\documentclass[12pt]{article}
\usepackage{amsmath,amscd,amsthm,amssymb}
\usepackage[russian]{babel}

\usepackage{a4wide}
\begin{document}
УДК 517.956

\begin{center}
\bf О НЕКЛАССИЧЕСКОЙ ТРАКТОВКЕ ЧЕТЫРЕХМЕРНОЙ ЗАДАЧИ ГУРСА ДЛЯ ОДНОГО ГИПЕРБОЛИЧЕСКОГО УРАВНЕНИЯ

\

И.Г. Мамедов

{\it (AZ 1141, г. Баку, ул. Б. Вагабзаде, 9, Институт Кибернетики им. А.И.Гусейнова НАН Азербайджана, www.cyber.az, e-mail: ilgar-mammadov@rambler.ru)}
\end{center}

\begin{abstract}
В данной статье выявлен гомеоморфизм между определенными парами банаховых пространств при исследовании четырехмерной задачи Гурса для одного дифференциального уравнения со старшей частной производной шестого порядка $D_{1} D_{2} D_{3}^{2} D_{4}^{2} u(x)$ с разрывными коэффициентами ($L_p-$коэффициентами) на основе сведения этой задачи к эквивалентному интегральному уравнению.

Библ. 24.

\

{\bf Ключевые слова:} гиперболическое уравнение, четырехмерная задача Гурса, уравнения с разрывными коэффициентами.
\end{abstract}

К настоящему времени усилиями многих математиков разнообразные классы трехмерных, четырехмерных, а также многомерных локальных и нелокальных начально-краевых задач для уравнений со старшей частной производной развивались в работах  [1]-[9]. Это связано с их появлением в различных задачах прикладного характера [10].\\

\centerline{\textbf{1. ПОСТАНОВКА ЗАДАЧИ} }

В работе обосновывается четырехмерная задача Гурса с неклассическими условиями для одного гиперболического уравнения.

Рассмотрим уравнение:
\[(V_{1,\, 1,\, 2,\, 2} u)(x)\equiv D_{1} D_{2} D_{3}^{2} D_{4}^{2} u(x)+a_{1,\, 0,\, 2,\, 2} (x)D_{1} D_{3}^{2} D_{4}^{2} u(x)+a_{0,\, 1,\, 2,\, 2} (x)D_{2} D_{3}^{2} D_{4}^{2} u(x)+
\]
\[
+a_{1,\, 1,\, 1,\, 2} (x)D_{1} D_{2} D_{3} D_{4}^{2} u(x)+a_{1,\, 1,\, 2,\, 1} (x)D_{1} D_{2} D_{3}^{2} D_{4} u(x)+
\]
\begin{equation} \label{GrindEQ__1_}
+\mathop{\sum }\limits_{\begin{array}{l} {i_{1} +i_{2} +i_{3} +i_{4} <5} \\ {i_{\xi } =\overline{0,\, 1},\, \; \xi =\overline{1,\, 2};} \\ {i_\eta =\overline{0,\, 2},\, \; \eta =\overline{3,\, 4}} \end{array}} a_{i_{1} ,\, i_{2} ,\, i_{3} ,\, i_{4} } (x)D_{1}^{i_{1} } D_{2}^{i_{2} } D_{3}^{i_{3} } D_{4}^{i_{4} } u(x)=\varphi _{1,\, 1,\, 2,\, 2} (x)\in L_{p} (G),
\end{equation}
здесь  $u(x)=u\left(x_{1} ,\, x_{2} ,\, x_{3} ,\, x_{4} \right)$  искомая функция, определенная на  $G;$  $a_{i_{1} ,\, i_{2} ,\, i_{3} ,\, i_{4} }=$ \linebreak$=a_{i_{1} ,\, i_{2} ,\, i_{3} ,\, i_{4} } (x)$ заданные измеримые функции на $G=G_{1} \times G_{2} \times G_{3} \times G_{4} $, где $G_{\xi } =(0,\, h_{\xi } )$, $\xi =\overline{1,\, 4};\, \varphi _{1,1,2,2} (x)$ заданная измеримая функция на $G;\, D_{\xi } =\dfrac{\partial }{\partial x_{\xi } } -$оператор обобщенного дифференцирования в смысле С.Л.Соболева.

Уравнение (1) является гиперболическим уравнением, которое обладает четырьмя действительными характеристиками $x_{1}=const,\, x_{2}=const,\, x_{3}=const,\, x_{4} =const,$ первая и вторая, из которых простая, а третая  и четвертая--двухкратная. Поэтому уравнение (1) в некотором смысле можно рассматривать также как псевдопараболическое уравнение [11]-[12]. Уравнения подобного вида возникают при описании многих реальных процессов, происходящих в природе и технике. Подобные ситуации имеют место при  изучении процессов распространения тепла [13], влагопереноса в почвогрунтах [14], фильтрации жидкости в пористых средах [15], в задачах математической биологии [16], а также в теории оптимальных процессов [17].

В этой работе уравнение (1) исследовано в общем случае, когда коэффициенты  $a_{i_{1} ,\, i_{2} ,\, i_{3} ,\, i_{4} } (x)$  являются негладкими функциями, удовлетворяющими лишь следующим условиям:
\[
a_{0,\, 0,\, i_{3} ,\, i_{4} } (x)\in L_{p} (G),\, \; \; a_{1,\, 0,\, i_{3} ,\, i_{4} } (x)\in L_{\infty ,\, p,\, p,\, p}^{x_{1} ,\, x_{2} ,\, x_{3} ,\, x_{4} } (G),a_{0,\, 1,\, i_{3} ,\, i_{4} } (x)\in L_{p,\infty ,\, p,\, p}^{x_{1} ,\, x_{2} ,\, x_{3} ,\, x_{4} } (G),
 \]\[
 a_{0,\, 0,\, 2,\, i_{4} } (x)\in L_{p,\, p,\infty ,\, p}^{x_{1} ,\, x_{2} ,\, x_{3} ,\, x_{4} } (G),a_{0,\, 0,\, i_{3} ,\, 2} (x)\in L_{p,\, p,\, p,\infty }^{x_{1} ,\, x_{2} ,\, x_{3} ,\, x_{4} } (G),\, \, a_{1,\, 1,\, i_{3} ,\, i_{4} } (x)\in L_{\infty ,\infty ,\, p,\, p}^{x_{1} ,\, x_{2} ,\, x_{3} ,\, x_{4} } (G),
 \]
 \[
 a_{1,\, 0,\, 2,\, i_{4} } (x)\in L_{\infty ,\, p,\infty ,\, p}^{x_{1} ,\, x_{2} ,\, x_{3} ,\, x_{4} } (G),\; a_{1,\, 0,\, i_{3} ,\, 2} (x)\in L_{\infty ,\, p,\, p,\infty }^{x_{1} ,\, x_{2} ,\, x_{3} ,\, x_{4} } (G),
a_{0,\, 1,\, 2,\, i_{4} } (x)\in L_{p,\infty ,\infty ,\, p}^{x_{1} ,\, x_{2} ,\, x_{3} ,\, x_{4} } (G),
\]
\[
a_{0,\, 1,\, i_{3} ,\, 2} (x)\in L_{p,\infty ,\, p,\infty }^{x_{1} ,\, x_{2} ,\, x_{3} ,\, x_{4} } (G), a_{0,\, 0,\, 2,\, 2} (x)\in L_{p,\, p,\infty ,\infty }^{x_{1} ,\, x_{2} ,\, x_{3} ,\, x_{4} } (G),\, a_{1,\, 1,\, 2,\, i_{4} } (x)\in L_{\infty ,\infty ,\infty ,\, p}^{x_{1} ,\, x_{2} ,\, x_{3} ,\, x_{4} } (G),
\]
\[
a_{1,\, 1,\, i_{3} ,\, 2} (x)\in L_{\infty ,\infty ,\, p,\infty }^{x_{1} ,\, x_{2} ,\, x_{3} ,\, x_{4} } (G),\, \; a_{0,\, 1,\, 2,\, 2} (x)\in L_{p,\infty ,\infty ,\infty }^{x_{1} ,\, x_{2} ,\, x_{3} ,\, x_{4} } (G), a_{1,\, 0,\, 1,\, 1} (x)\in L_{\infty ,\, p,\infty ,\infty }^{x_{1} ,\, x_{2} ,\, x_{3} ,\, x_{4} } (G),
\]
где $i_{3} =\overline{0,\, 1},\, \; i_{4} =\overline{0,\, 1}.$

При этом важным принципиальным моментом является то, что рассматриваемое уравнение обладает, разрывными коэффициентами которые удовлетворяют только некоторым условиям типа $p$-интегрируемости и ограниченности т.е. рассмотренный псевдопараболический дифференциальный оператор не имеет традиционного сопряженного оператора. Иначе говоря, функция Римана для такого уравнения не может быть исследована классическим методом характеристик.

При этих условиях решение $u(x)$ уравнения (1) будем искать в пространстве\linebreak С.Л.Соболева
\[
W_{p}^{(1,1,2,2)} (G)\equiv \left\{u(x): D_{1}^{i_{1} } D_{2}^{i_{2} } D_{3}^{i_{3} } D_{4}^{i_{4} } u(x)\in L_{p} (G),\; i_{\xi } =\overline{0, 1}, \xi =\overline{1, 2}; i_{\eta } =\overline{0, 2},  \eta =\overline{3,4}\right\},
\]
где $1\le p\le \infty $. Норму в анизотропном пространстве $W_{p}^{(1,1,2,2)} (G)$ будем определять равенством
\[
\left\| u(x)\right\| _{W_{p}^{(1,\, 1,\, 2,\, 2)} (G)} =\mathop{\sum }\limits_{\begin{array}{l} {i_{\xi } =0} \\ {\xi =\overline{1,\, 2}} \end{array}}^{1} \; \mathop{\sum }\limits_{\begin{array}{l} {i_{\eta } =0} \\ {\eta =\overline{3,\, 4}} \end{array}}^{2} \left\| D_{1}^{i_{1} } D_{2}^{i_{2} } D_{3}^{i_{3} } D_{4}^{i_{4} } u(x)\right\| _{L_{p} (G)} .
\]

Для уравнения (1) условия Гурса классического вида можно задать в виде
\begin{equation} \label{GrindEQ__2_}
\left\{\begin{array}{l} {u\left. \left(x_{1} ,\, x_{2} ,\, x_{3} ,\, x_{4} \right)\right|_{x_{1} =0} =F\left(x_{2} ,\, x_{3} ,\, x_{4} \right);\; \; \; u\left. \left(x_{1} ,\, x_{2} ,\, x_{3} ,\, x_{4} \right)\right|_{x_{2} =0} =g\left(x_{1} ,\, x_{3} ,\, x_{4} \right);} \\ {u\left. \left(x_{1} ,\, x_{2} ,\, x_{3} ,\, x_{4} \right)\right|_{x_{3} =0} =\psi \left(x_{1} ,\, x_{2} ,\, x_{4} \right);\; \; \; \left. \frac{\partial u\left(x_{1} ,\, x_{2} ,\, x_{3} ,\, x_{4} \right)}{\partial x_{3} } \right|_{_{x_{3} =0} } =\Phi \left(x_{1} ,\, x_{2} ,\, x_{4} \right);} \\ {u\left. \left(x_{1} ,\, x_{2} ,\, x_{3} ,\, x_{4} \right)\right|_{x_{4} =0} =T\left(x_{1} ,\, x_{2} ,\, x_{3} \right);\; \; \; \left. \frac{\partial u\left(x_{1} ,\, x_{2} ,\, x_{3} ,\, x_{4} \right)}{\partial x_{4} } \right|_{_{x_{4} =0} } =S\left(x_{1} ,\, x_{2} ,\, x_{3} \right).} \end{array}\right.
\end{equation}
где $F\left(x_{2} , x_{3} , x_{4} \right),\, g\left(x_{1} , x_{3} , x_{4} \right), \psi \left(x_{1} , x_{2} , x_{4} \right), \Phi \left(x_{1} , x_{2} , x_{4} \right), T\left(x_{1} , x_{2} , x_{3} \right), S\left(x_{1} , x_{2} , x_{3} \right)$, заданные изме-римые  функции на $G.$ Очевидно, что в случае условий (2) функций $F,\, g,\, \psi ,\, \Phi ,\, T,\, S$ кроме условий
\[
F\left(x_{2} ,\, x_{3} ,\, x_{4} \right)\in W_{p}^{(1,\, 2,\, 2)} \left(G_{2} \times G_{3} \times G_{4} \right),   g\left(x_{1} ,\, x_{3} ,\, x_{4} \right)\in W_{p}^{(1,\, 2,\, 2)} \left(G_{1} \times G_{3} \times G_{4} \right),
\]
\[
\psi \left(x_{1} ,\, x_{2} ,\, x_{4} \right)\in W_{p}^{(1,\, 1,\, 2)} \left(G_{1} \times G_{2} \times G_{4} \right),   \Phi \left(x_{1} ,\, x_{2} ,\, x_{4} \right)\in W_{p}^{(1,\, 1,\, 2)} \left(G_{1} \times G_{2} \times G_{4} \right),
\]
\[
T\left(x_{1} ,\, x_{2} ,\, x_{3} \right)\in W_{p}^{(1,\, 1,\, 2)} \left(G_{1} \times G_{2} \times G_{3} \right),  S\left(x_{1} ,\, x_{2} ,  x_{3} \right)\in W_{p}^{(1,  1,  2)} \left(G_{1} \times G_{2} \times G_{3} \right),
\]
должны удовлетворять также следующим условиям:
\begin{equation} \label{GrindEQ__3_}
\left\{\!\!\!\begin{array}{l} {F\left(0, x_{3}, x_{4} \right)\!=\!g\left(0, x_{3},  x_{4} \right);    F\left(x_{2},  0,  x_{4} \right)\!=\!\psi \left(0,  x_{2},  x_{4} \right);    g_{x_{4} } (x_{1},x_{3} ,0)\!=\!S\left(x_{1},0,x_{3} \right);} \\
{F\left(x_{2}, \! x_{3},\!  0\right)\!=\!T\left(0, x_{2}, x_{3} \right)\!;\! F_{x_{3}} \left(x_{2},  0,  x_{4} \right)\!=\!\Phi \left(0,  x_{2},  x_{4}\right)\!;\! F_{x_{4} } \left(x_{2},  x_{3},0\right)\!=\!S\left(0,  x_{2},  x_{3} \right);} \\
{g\left(x_{1},  x_{3},0\right)\!=\!T\left(x_{1},  0,  x_{3} \right);  g\left(x_{1},0,x_{4} \right)\!=\!\psi \left(x_{1},0,x_{4} \right);    g_{x_{3} } (x_{1},0,x_{4} )\!=\!\Phi \left(x_{1},0,x_{4} \right);} \\
{\psi (x_{1},x_{2} ,0)\!=\!T\left(x_{1} ,x_{2},0\right)\!; \!   \psi _{x_{4} } (x_{1},  x_{2} ,0)\!=\!S\left(x_{1},  x_{2}, 0\right)\!;\! \Phi _{x_{4} } \left(x_{1} ,x_{2},0\right)\!=\!S_{x_{3} } (x_{1},  x_{2},0)}
\end{array}
\right.
\end{equation}
которые являются условиями согласования.

Наличие условий согласования (3) в постановке задачи (1), (2) означает, что условиями (2) задана также некоторая излишняя информация о решении этой задачи. Поэтому возникает вопрос о нахождении краевых условий, которые не содержат излишней информации о решении и не требуют выполнения некоторых дополнительных условий типа согласования. В связи с этим рассмотрим следующие неклассические начально-краевые условия:
\[
{V_{0,0,i_{3} ,i_{4} } u\equiv D_{3}^{i_{3} } D_{4}^{i_{4} } u\left(0,0,0,0\right)=\varphi _{0,0,i_{3} ,i_{4} } \in R,              i_{\nu } =\overline{0,1},           \     \nu =\overline{3,4};}
\]
\[
{\left(V_{1,0,i_{3} ,i_{4} } u\right)\left(x_{1} \right)\equiv D_{1} D_{3}^{i_{3} } D_{4}^{i_{4} } u\left(x_{1} ,0,0,0\right)=\varphi _{1,0,i_{3} ,i_{4} } \left(x_{1} \right)\in L_{p} \left(G_{1} \right),              i_{\nu } =\overline{0,1},
        \nu =\overline{3,4};}
\]
\[
\left(V_{0,1,i_{3} ,i_{4} } u\right)\left(x_{2} \right)\equiv D_{2} D_{3}^{i_{3} } D_{4}^{i_{4} } u\left(0,x_{2} ,0,0\right)=\varphi _{0,1,i_{3} ,i_{4} } \left(x_{2} \right)\in L_{p} \left(G_{2} \right),\ i_{\nu } =\overline{0,1}, \ \nu =\overline{3,4};
\]
\[
\left(V_{0,0,2,i_{4} } u\right)\left(x_{3} \right)\equiv D_{3}^{2} D_{4}^{i_{4} } u\left(0,0,x_{3} ,0\right)=\varphi _{0,0,2,i_{4} } \left(x_{3} \right)\in L_{p} \left(G_{3} \right),        i_{4} =\overline{0,1};
\]
\[
{\left(V_{0,0,i_{3} ,2} u\right)\left(x_{4} \right)\equiv D_{3}^{i_{3} } D_{4}^{2} u\left(0,0,0,x_{4} \right)=\varphi _{0,0,i_{3} ,2} \left(x_{4} \right)\in L_{p} \left(G_{4} \right),        i_{3} =\overline{0,1};}
\]
\[ \left(V_{1,1,i_{3} ,i_{4} } u\right)\left(x_{1} ,x_{2} \right)\equiv D_{1} D_{2} D_{3}^{i_{3} } D_{4}^{i_{4} } u\left(x_{1} ,x_{2} ,0,0\right)=
\]\[
=\varphi _{1,1,i_{3} ,i_{4} } \left(x_{1} ,x_{2} \right)\in L_{p} \left(G_{1} \times G_{2} \right),        i_{\nu } =\overline{0,1},      \nu =\overline{3,4};
\]
\[ {\left(V_{1,0,2,i_{4} } u\right)\left(x_{1} ,x_{3} \right)\equiv D_{1} D_{3}^{2} D_{4}^{i_{4} } u\left(x_{1} ,0,x_{3} ,0\right)=\varphi _{1,0,2,i_{4} } \left(x_{1} ,x_{3} \right)\in L_{p} \left(G_{1} \times G_{3} \right),              i_{4} =\overline{0,1};}
\]
\[
 {\left(V_{1,0,i_{3} ,2} u\right)\left(x_{1} ,x_{4} \right)\equiv D_{1} D_{3}^{i_{3} } D_{4}^{2} u\left(x_{1} ,0,0,x_{4} \right)=\varphi _{1,0,i_{3} ,2} \left(x_{1} ,x_{4} \right)\in L_{p} \left(G_{1} \times G_{4} \right),              i_{3} =\overline{0,1};}
 \]
 \[
  {\left(V_{0,1,2,i_{4} } u\right)\left(x_{2} ,x_{3} \right)\equiv D_{2} D_{3}^{2} D_{4}^{i_{4} } u\left(0,x_{2} ,x_{3} ,0\right)=\varphi _{0,1,2,i_{4} } \left(x_{2} ,x_{3} \right)\in L_{p} \left(G_{2} \times G_{3} \right),              i_{4} =\overline{0,1};}
  \]
\[ {\left(V_{0,1,i_{3} ,2} u\right)\left(x_{2} ,x_{4} \right)\equiv D_{2} D_{3}^{i_{3} } D_{4}^{2} u\left(0,x_{2} ,0,x_{4} \right)=\varphi _{0,1,i_{3} ,2} \left(x_{2} ,x_{4} \right)\in L_{p} \left(G_{2} \times G_{4} \right),              i_{3} =\overline{0,1};}
\]
\[ {\left(V_{0,0,2,2} u\right)\left(x_{3} ,x_{4} \right)\equiv D_{3}^{2} D_{4}^{2} u\left(0,0,x_{3} ,x_{4} \right)=\varphi _{0,0,2,2} \left(x_{3} ,x_{4} \right)\in L_{p} \left(G_{3} \times G_{4} \right);}
\]
\[
\left(V_{1,1,2,i_{4} } u\right)\left(x_{1} ,x_{2} ,x_{3} \right)\equiv D_{1} D_{2} D_{3}^{2} D_{4}^{i_{4} } u\left(x_{1} ,x_{2} ,x_{3} ,0\right)=
\]\[
=\varphi _{1,1,2,i_{4} } \left(x_{1} ,x_{2} ,x_{3} \right)\in L_{p} \left(G_{1} \times G_{2} \times G_{3} \right),    i_{4} =\overline{0,1};
\]
\[
\left(V_{1,1,i_{3} ,2} u\right)\left(x_{1} ,x_{2} ,x_{4} \right)\equiv D_{1} D_{2} D_{3}^{i_{3} } D_{4}^{2} u\left(x_{1} ,x_{2} ,0,x_{4} \right)=
\]\[
=\varphi _{1,1,i_{3} ,2} \left(x_{1} ,x_{2} ,x_{4} \right)\in L_{p} \left(G_{1} \times G_{2} \times G_{4} \right),      i_{3} =\overline{0,1};
\]
\[
\left(V_{0,1,2,2} u\right)\left(x_{2} ,x_{3} ,x_{4} \right)\equiv D_{2} D_{3}^{2} D_{4}^{2} u\left(0,x_{2} ,x_{3} ,x_{4} \right)=\varphi _{0,1,2,2} \left(x_{2} ,x_{3} ,x_{4} \right)\in L_{p} \left(G_{2} \times G_{3} \times G_{4} \right);
\]
\begin{equation} \label{GrindEQ__4_}
\left(V_{1,0,2,2} u\right)\left(x_{1} ,x_{3} ,x_{4} \right)\equiv D_{1} D_{3}^{2} D_{4}^{2} u\left(x_{1} ,0,x_{3} ,x_{4} \right)=
$$$$
=\varphi _{1,0,2,2} \left(x_{1} ,x_{3} ,x_{4} \right)\in L_{p} \left(G_{1} \times G_{3} \times G_{4} \right);
\end{equation}

\,Если функция $u\in W_{p}^{(1,  1,2,  2)} (G)$ является решением четырехмерной задачи Гурса классического вида (1), (2), то она является также решением задачи (1), (4) для  $\varphi _{i_{1} ,  i_{2} ,  i_{3} ,  i_{4} } $ определяемых следующими равенствами
\[
\varphi _{0,0,0,0} =F\left(0,0,0\right)=g\left(0,0,0\right)=\psi \left(0,0,0\right)=T\left(0,0,0\right);  \varphi _{0,0,1,0}=
\]\[
 =\Phi \left(0,0,0\right)=g_{x_{3} } \left(0,0,0\right)=F_{x_{3} } \left(0,0,0\right);
\]
\[
\varphi _{0,0,0,1} =S\left(0,0,0\right)=F_{x_{4} } \left(0,0,0\right)=\psi _{x_{4} } \left(0,0,0\right); \varphi _{0,0,1,1} =S_{x_{3} } \left(0,0,0\right)=\Phi _{x_{4} } \left(0,0,0\right);
\]
\[
\varphi _{1,0,0,0} \left(x_{1} \right)=g_{x_{1} } \left(x_{1} ,0,0\right)=\psi _{x_{1} } \left(x_{1} ,0,0\right)
=T_{x_{1} } \left(x_{1} ,0,0\right);\varphi _{1,0,1,0} \left(x_{1} \right)
=g_{x_{1} x_{3} } \left(x_{1} ,0,0\right)=
\]
\[
=\Phi _{x_{1} } \left(x_{1} ,0,0\right)
=T_{x_{1} x_{3} } \left(x_{1} ,0,0\right);\varphi _{1,0,0,1} \left(x_{1} \right)=g_{x_{1} x_{4} } \left(x_{1} ,0,0\right)=\psi _{x_{1} x_{4} } \left(x_{1} ,0,0\right)
=
\]\[
=S_{x_{1} } \left(x_{1} ,0,0\right); \
\varphi _{1,0,1,1} \left(x_{1} \right)=g_{x_{1} x_{3} x_{4} } \left(x_{1} ,0,0\right)
=\Phi _{x_{1} x_{4} } \left(x_{1} ,0,0\right);   \varphi _{0,1,0,0} \left(x_{2} \right)=
\]
\[
=F_{x_{2} } \left(x_{2} ,0,0\right)=\psi _{x_{2} } \left(0,x_{2} ,0\right)=T_{x_{2} } \left(0,x_{2} ,0\right);\varphi _{0,1,1,0} \left(x_{2} \right)=F_{x_{2} x_{3} } \left(x_{2} ,0,0\right)=
\]\[
=T_{x_{2} x_{3} } \left(0,x_{2} ,0\right)=\Phi _{x_{2} } \left(0,x_{2} ,0\right);
\varphi _{0,1,0,1} \left(x_{2} \right)=F_{x_{2} x_{4} } \left(x_{2} ,0,0\right)=\psi _{x_{2} x_{4} } \left(0,x_{2} ,0\right)=
\]\[
=S_{x_{2} } \left(0,x_{2} ,0\right); \varphi _{0,1,1,1} \left(x_{2} \right)=F_{x_{2} x_{3} x_{4} } \left(x_{2} ,0,0\right)=S_{x_{2} x_{3} } \left(0,x_{2} ,0\right);
\]
\[
\varphi _{0,0,2,0} \left(x_{3} \right)=F_{x_{3} x_{3} } \left(0,x_{3} ,0\right)=g_{x_{3} x_{3} } \left(0,x_{3} ,0\right)=T_{x_{3} x_{3} } \left(0,0,x_{3} \right);
\]\[
\varphi _{0,0,2,1} \left(x_{3} \right)=F_{x_{3} x_{3} x_{4} } \left(0,x_{3} ,0\right)=g_{x_{3} x_{3} x_{4} } \left(0,x_{3} ,0\right)=S_{x_{3} x_{3} } \left(0,0,x_{3} \right);
\]
\[
\varphi _{0,0,0,2} \left(x_{4} \right)=F_{x_{4} x_{4} } \left(0,0,x_{4} \right)=g_{x_{4} x_{4} } \left(0,0,x_{4} \right)=\psi _{x_{4} x_{4} } \left(0,0,x_{4} \right);
\]
\[
\varphi _{0,0,1,2} \left(x_{4} \right)=F_{x_{3} x_{4} x_{4} } \left(0,0,x_{4} \right)=g_{x_{3} x_{4} x_{4} } \left(0,0,x_{4} \right)=\Phi _{x_{4} x_{4} } \left(0,0,x_{4} \right);
\]
\[
\varphi _{1,1,0,0} \left(x_{1} ,x_{2} \right)=\psi _{x_{1} x_{2} } \left(x_{1} ,x_{2} ,0\right)=T_{x_{1} x_{2} } \left(x_{1} ,x_{2} ,0\right);
\]\[
\varphi _{1,1,1,0} \left(x_{1} ,x_{2} \right)=T_{x_{1} x_{2} x_{3} } \left(x_{1} ,x_{2} ,0\right)=\Phi _{x_{1} x_{2} } \left(x_{1} ,x_{2} ,0\right);
\]
\[\varphi _{1,1,0,1} \left(x_{1} ,x_{2} \right)=\psi _{x_{1} x_{2} x_{4} } \left(x_{1} ,x_{2} ,0\right)=S_{x_{1} x_{2} } \left(x_{1} ,x_{2} ,0\right);
 \]\[
 \varphi _{1,1,1,1} \left(x_{1} ,x_{2} \right)=\Phi _{x_{1} x_{2} x_{4} } \left(x_{1} ,x_{2} ,0\right)=S_{x_{1} x_{2} x_{3} } \left(x_{1} ,x_{2} ,0\right);
\]
\[\varphi _{1,0,2,0} \left(x_{1} ,x_{3} \right)=g_{x_{1} x_{3} x_{3} } \left(x_{1} ,x_{3} ,0\right)=T_{x_{1} x_{3} x_{3} } \left(x_{1} ,0,x_{3} \right);
 \]\[
 \varphi _{1,0,2,1} \left(x_{1} ,x_{3} \right)=g_{x_{1} x_{3} x_{3} x_{4} } \left(x_{1} ,x_{3} ,0\right)=S_{x_{1} x_{3} x_{3} } \left(x_{1} ,0,x_{3} \right);
\]
\[\varphi _{1,0,0,2} \left(x_{1} ,x_{4} \right)=g_{x_{1} x_{4} x_{4} } \left(x_{1} ,0,x_{4} \right)=\psi _{x_{1} x_{4} x_{4} } \left(x_{1} ,0,x_{4} \right);
 \]
 \[\varphi _{1,0,1,2} \left(x_{1} ,x_{4} \right)=g_{x_{1} x_{3} x_{4} x_{4} } \left(x_{1} ,0,x_{4} \right)=\Phi _{x_{1} x_{4} x_{4} } \left(x_{1} ,0,x_{4} \right);
\]
\[
\varphi _{0,1,2,0} \left(x_{2} ,x_{3} \right)=F_{x_{2} x_{3} x_{3} } \left(x_{2} ,x_{3} ,0\right)=T_{x_{2} x_{3} x_{3} } \left(0,x_{2} ,x_{3} \right);
\]\[
\varphi _{0,1,2,1} \left(x_{2} ,x_{3} \right)=F_{x_{2} x_{3} x_{3} x_{4} } \left(x_{2} ,x_{3} ,0\right)=S_{x_{2} x_{3} x_{3} } \left(0,x_{2} ,x_{3} \right);
\]
\[
\varphi _{0,1,0,2} \left(x_{2} ,x_{4} \right)=F_{x_{2} x_{4} x_{4} } \left(x_{2} ,0,x_{4} \right)=\psi _{x_{2} x_{4} x_{4} } \left(0,x_{2} ,x_{4} \right);
\]\[
\varphi _{0,1,1,2} \left(x_{2} ,x_{4} \right)=F_{x_{2} x_{3} x_{4} x_{4} } \left(x_{2} ,0,x_{4} \right)=\Phi _{x_{2} x_{4} x_{4} } \left(0,x_{2} ,x_{4} \right);
\]
\[
\varphi _{0,0,2,2} \left(x_{3} ,x_{4} \right)=F_{x_{3} x_{3} x_{4} x_{4} } \left(0,x_{3} ,x_{4} \right)=g_{x_{3} x_{3} x_{4} x_{4} } \left(0,x_{3} ,x_{4} \right); \varphi _{1,1,2,0} \left(x_{1} ,x_{2} ,x_{3} \right)=
\]\[
=T_{x_{1} x_{2} x_{3} x_{3} } \left(x_{1} ,x_{2} ,x_{3} \right); \varphi _{1,1,2,1} \left(x_{1} ,x_{2} ,x_{3} \right)=S_{x_{1} x_{2} x_{3} x_{3} } \left(x_{1} ,x_{2} ,x_{3} \right);  \varphi _{1,1,0,2} \left(x_{1} ,x_{2} ,x_{4} \right)=
\]\[
=\psi _{x_{1} x_{2} x_{4} x_{4} } \left(x_{1} ,x_{2} ,x_{4} \right);
\
\varphi _{1,1,1,2} \left(x_{1} ,x_{2} ,x_{4} \right)=\Phi _{x_{1} x_{2} x_{4} x_{4} } \left(x_{1} ,x_{2} ,x_{4} \right); \]
\[ \varphi _{0,1,2,2} \left(x_{2} ,x_{3} ,x_{4} \right)=F_{x_{2} x_{3} x_{3} x_{4} x_{4} } \left(x_{2} ,x_{3} ,x_{4} \right);
 \varphi _{1,0,2,2} \left(x_{1} ,x_{3} ,x_{4} \right)=g_{x_{1} x_{3} x_{3} x_{4} x_{4} } \left(x_{1} ,x_{3} ,x_{4} \right).
\]

Легко доказать, что верно и обратное. Иначе говоря, если функция $u\in W_{p}^{\left(1,1,2,2\right)} \left(G\right)$ является решением задачи (1), (4), то она является также решением задачи (1), (2), для следующих функций $F,      g,      \psi ,      \Phi ,      T,      S  $:
\[
F\left(x_{2} ,x_{3} ,x_{4} \right)=\sum _{i_{3} =0}^{1}  \sum _{i_{4} =0}^{1}x_{3}^{i_{3} } x_{4}^{i_{4} } \varphi _{0,0,i_{3} ,i_{4} }   +\sum _{i_{3} =0}^{1}  \sum _{i_{4} =0}^{1}x_{3}^{i_{3} } x_{4}^{i_{4} } \int _{0}^{x_{2} }\varphi _{0,1,i_{3} ,i_{4} } \left(\xi _{2} \right)d\xi _{2} +
   \]
   \[
   +\sum _{i_{4} =0}^{1}  x_{4}^{i_{4} } \int _{0}^{x_{3} }\left(x_{3} -\xi _{3} \right)\varphi _{0,0,2,i_{4} } \left(\xi _{3} \right)d\xi _{3}
+\sum _{i_{3} =0}^{1}  x_{3}^{i_{3} } \int _{0}^{x_{4} }\left(x_{4} -\xi _{4} \right)\varphi _{0,0,i_{3} ,2} \left(\xi _{4} \right)d\xi _{4} +
  \]\[
  +\sum _{i_{4} =0}^{1}  x_{4}^{i_{4} } \int _{0}^{x_{2} }\int _{0}^{x_{3} }\left(x_{3} -\xi _{3} \right)\varphi _{0,1,2,i_{4} } \left(\xi _{2} ,\xi _{3} \right)d\xi _{2} d\xi _{3}    +
\]\[
+\int _{0}^{x_{3} }\int _{0}^{x_{4} }\left(x_{3} -\xi _{3} \right)  \left(x_{4} -\xi _{4} \right)  \varphi _{0,0,2,2} \left(\xi _{3} ,\xi _{4} \right)d\xi _{3} d\xi _{4} +
  \]\[
  +\sum _{i_{3} =0}^{1}  x_{3}^{i_{3} } \int _{0}^{x_{2} }\int _{0}^{x_{4} }\left(x_{4} -\xi _{4} \right)\varphi _{0,1,i_{3} ,2} \left(\xi _{2} ,\xi _{4} \right)d\xi _{2} d\xi _{4}    +
\]
\[
+\int _{0}^{x_{2} }\int _{0}^{x_{3} }\int _{0}^{x_{4} }\left(x_{3} -\xi _{3} \right)  \left(x_{4} -\xi _{4} \right)  \varphi _{0,1,2,2} \left(\xi _{2} ,\xi _{3} ,\xi _{4} \right)d\xi _{2} d\xi _{3} d\xi _{4}    ;
\]
\[
g\left(x_{1} ,x_{3} ,x_{4} \right)=\sum _{i_{3} =0}^{1}  \sum _{i_{4} =0}^{1}x_{3}^{i_{3} } x_{4}^{i_{4} } \varphi _{0,0,i_{3} ,i_{4} }   +\sum _{i_{3} =0}^{1}  \sum _{i_{4} =0}^{1}x_{3}^{i_{3} } x_{4}^{i_{4} } \int _{0}^{x_{1} }\varphi _{1,0,i_{3} ,i_{4} } \left(\xi _{1} \right)d\xi _{1}
  +
  \]\[
  +\sum _{i_{4} =0}^{1}  x_{4}^{i_{4} } \int _{0}^{x_{3} }\left(x_{3} -\xi _{3} \right)\varphi _{0,0,2,i_{4} } \left(\xi _{3} \right)d\xi _{3}
+\sum _{i_{3} =0}^{1}  x_{3}^{i_{3} } \int _{0}^{x_{4} }\left(x_{4} -\xi _{4} \right)\varphi _{0,0,i_{3} ,2} \left(\xi _{4} \right)d\xi _{4} +
  \]
  \[+\sum _{i_{4} =0}^{1}  x_{4}^{i_{4} } \int _{0}^{x_{1} }\int _{0}^{x_{3} }\left(x_{3} -\xi _{3} \right)\varphi _{1,0,2,i_{4} } \left(\xi _{1} ,\xi _{3} \right)d\xi _{1} d\xi _{3}    +
\]
\[
+\int _{0}^{x_{3} }\int _{0}^{x_{4} }\left(x_{3} -\xi _{3} \right)  \left(x_{4} -\xi _{4} \right)  \varphi _{0,0,2,2} \left(\xi _{3} ,\xi _{4} \right)d\xi _{3} d\xi _{4} +
  \]\[
  +\sum _{i_{3} =0}^{1}  x_{3}^{i_{3} } \int _{0}^{x_{1} }\int _{0}^{x_{4} }\left(x_{4} -\xi _{4} \right)\varphi _{1,0,i_{3} ,2} \left(\xi _{1} ,\xi _{4} \right)d\xi _{1} d\xi _{4}    +
\]
\[
+\int _{0}^{x_{1} }\int _{0}^{x_{3} }\int _{0}^{x_{4} }\left(x_{3} -\xi _{3} \right)  \left(x_{4} -\xi _{4} \right)  \varphi _{1,0,2,2} \left(\xi _{1} ,\xi _{3} ,\xi _{4} \right)d\xi _{1} d\xi _{3} d\xi _{4}    ;
\]
$\psi \left(x_{1} ,x_{2} ,x_{4} \right)=M_{0} \left(x_{1} ,x_{2} ,x_{4} \right),\; \; \Phi \left(x_{1} ,x_{2} ,x_{4} \right)=M_{1} \left(x_{1} ,x_{2} ,x_{4} \right)$, где
\[
M_{k} \left(x_{1} ,x_{2} ,x_{4} \right)=\varphi _{0,0,k,0} +x_{4} \varphi _{0,0,k,1} +\int _{0}^{x_{1} }\varphi _{1,0,k,0} \left(\tau _{1} \right)d\tau _{1} +
x_{4} \int _{0}^{x_{1} }\varphi _{1,0,k,1} \left(\tau _{1} \right)d\tau _{1}  +
\]\[
 +\int _{0}^{x_{2} }\varphi _{0,1,k,0} \left(\tau _{2} \right)d\tau _{2}  +
x_{4} \int _{0}^{x_{2} }\varphi _{0,1,k,1} \left(\tau _{2} \right)d\tau _{2}  +\int _{0}^{x_{4} }\left(x_{4} -\tau _{4} \right)\varphi _{0,0,k,2} \left(\tau _{4} \right)d\tau _{4}+
  \]\[
  +\int _{0}^{x_{1} }\int _{0}^{x_{2} }\varphi _{1,1,k,0} \left(\tau _{1} ,\tau _{2} \right)d\tau _{1} d\tau _{2}   +x_{4} \int _{0}^{x_{1} }\int _{0}^{x_{2} }\varphi _{1,1,k,1} \left(\tau _{1} ,\tau _{2} \right)d\tau _{1} d\tau _{2}   +
\]
\[
+\int _{0}^{x_{2} }\int _{0}^{x_{4} }\left(x_{4} -\tau _{4} \right)\varphi _{0,1,k,2} \left(\tau _{2} ,\tau _{4} \right)d\tau _{2} d\tau _{4}   +\int _{0}^{x_{1} }\int _{0}^{x_{4} }\left(x_{4} -\tau _{4} \right)\varphi _{1,0,k,2} \left(\tau _{1} ,\tau _{4} \right)d\tau _{1} d\tau _{4}   +
\]
\[
+\int _{0}^{x_{1} }\int _{0}^{x_{2} }\int _{0}^{x_{4} }\left(x_{4} -\tau _{4} \right)\varphi _{1,1,k,2} \left(\tau _{1} ,\tau _{2} ,\tau _{4} \right)d\tau _{1} d\tau _{2} d\tau _{4}      ,\; \; \; k=\overline{0,1}    ;
\]
$T\left(x_{1} ,x_{2} ,x_{3} \right)=L_{0} \left(x_{1} ,x_{2} ,x_{3} \right)\; ,\; \; S\left(x_{1} ,x_{2} ,x_{3} \right)=L_{1} \left(x_{1} ,x_{2} ,x_{3} \right)$,  где
\[
L_{k} \left(x_{1} ,x_{2} ,x_{3} \right)=\varphi _{0,0,0,k} +x_{3} \varphi _{0,0,1,k} +\int _{0}^{x_{1} }\varphi _{1,0,0,k} \left(\eta _{1} \right)d\eta _{1} +
 \]\[
 +x_{3} \int _{0}^{x_{1} }\varphi _{1,0,1,k} \left(\eta _{1} \right)d\eta _{1}  +\int _{0}^{x_{2} }\varphi _{0,1,0,k} \left(\eta _{2} \right)d\eta _{2}  +
\]
\[
+x_{3} \int _{0}^{x_{2} }\varphi _{0,1,1,k} \left(\eta _{2} \right)d\eta _{2}  +\int _{0}^{x_{3} }\left(x_{3} -\eta _{3} \right)\varphi _{0,0,2,k} \left(\eta _{3} \right)d\eta _{3}  +\int _{0}^{x_{1} }\int _{0}^{x_{2} }\varphi _{1,1,0,k} \left(\eta _{1} ,\eta _{2} \right)d\eta _{1} d\eta _{2}   +
\]
\[
+x_{3} \int _{0}^{x_{1} }\int _{0}^{x_{2} }\varphi _{1,1,1,k} \left(\eta _{1} ,\eta _{2} \right)d\eta _{1} d\eta _{2}   +\int _{0}^{x_{2} }\int _{0}^{x_{3} }\left(x_{3} -\eta _{3} \right)\varphi _{0,1,2,k} \left(\eta _{2} ,\eta _{3} \right)d\eta _{2} d\eta _{3}   +
\]
\[
+\int _{0}^{x_{1} }\int _{0}^{x_{3} }\left(x_{3} -\eta _{3} \right)\varphi _{1,0,2,k} \left(\eta _{1} ,\eta _{3} \right)d\eta _{1} d\eta _{3} +
  \]\[
  +\int _{0}^{x_{1} }\int _{0}^{x_{2} }\int _{0}^{x_{3} }\left(x_{3} -\eta _{3} \right)\varphi _{1,1,2,k} \left(\eta _{1} ,\eta _{2} ,\eta _{3} \right)d\eta _{1} d\eta _{2} d\eta _{3}      ,\; \; \; k=\overline{0,1}.
\]

Итак, четырехмерной задачи Гурса классического вида (1), (2) и вида (1), (4) в общем случае эквивалентны. Однако четырехмерная задача Гурса (1), (4) по постановке более естественна, чем задача (1), (2). Это связано с тем, что в постановке задачи (1), (4) на правые части краевых условий никаких дополнительных условий типа согласования не требуется. Поэтому задачу (1), (4) можно рассматривать как задачу Гурса с неклассическими условиями.\\

\begin{center}
\bf 2. НЕКОТОРЫЕ СТРУКТУРНЫЕ СВОЙСТВА ПРОСТРАНСТВА С.Л.СОБОЛЕВА $W^{1,1,2,2}_p(G)$ И ОПЕРАТОРНЫЙ ВИД ЧЕТЫРЕХМЕРНОЙ ЗАДАЧИ ГУРСА
\end{center}

Задачу (1), (4) мы будем исследовать методом операторных уравнений и при этом будем следовать схемы работы [18]. Предварительно задачу (1), (4) запишем в виде операторного уравнения
\begin{equation} \label{GrindEQ__5_}
Vu=\varphi ,
\end{equation}
где $V$ есть векторный оператор, определяемый посредством равенства
\[
V=\left(V_{1,1,2,2} ,V_{0,0,0,0} ,V_{0,0,1,0} ,V_{0,0,0,1} ,V_{0,0,1,1} ,V_{1,0,0,0} ,V_{1,0,1,0} ,V_{1,0,0,1} ,V_{1,0,1,1} ,V_{0,1,0,0} ,V_{0,1,1,0} ,\right.
\]
\[
V_{0,1,0,1} ,V_{0,1,1,1} ,V_{0,0,2,0} ,V_{0,0,2,1} ,V_{0,0,0,2} ,V_{0,0,1,2} ,V_{1,1,0,0} ,V_{1,1,1,0} ,V_{1,1,0,1} ,V_{1,1,1,1} ,V_{1,0,2,0} ,V_{1,0,2,1}  ,
\]
\[
V_{1,0,0,2} ,V_{1,0,1,2},V_{0,1,2,0} ,V_{0,1,2,1} ,V_{0,1,0,2} ,V_{0,1,1,2} ,V_{0,0,2,2} ,V_{1,1,2,0} ,V_{1,1,2,1} ,V_{1,1,0,2} ,V_{1,1,1,2} ,
\]
$\left. V_{0,1,2,2} ,V_{1,0,2,2} \right):W_{p}^{\left(1,1,2,2\right)} \left(G\right)\to E_{p}^{\left(1,1,2,2\right)} $ а $\varphi $
есть заданный векторный элемент вида
\[
\varphi =\left(\varphi _{1,1,2,2} ,\varphi _{0,0,0,0} ,\varphi _{0,0,1,0} ,\varphi _{0,0,0,1} ,\varphi _{0,0,1,1} ,\varphi _{1,0,0,0} ,\varphi _{1,0,1,0} ,\varphi _{1,0,0,1} ,\varphi _{1,0,1,1} ,\varphi _{0,1,0,0} ,\right.
\]
\[
\varphi _{0,1,1,0} ,\varphi _{0,1,0,1} ,\varphi _{0,1,1,1} ,\varphi _{0,0,2,0} ,\varphi _{0,0,2,1} ,\varphi _{0,0,0,2} ,\varphi _{0,0,1,2} ,\varphi _{1,1,0,0} ,\varphi _{1,1,1,0} ,\varphi _{1,1,0,1} ,\varphi _{1,1,1,1} ,
\]
\[
\varphi _{1,0,2,0} ,\varphi _{1,0,2,1} ,\varphi _{1,0,0,2} ,\varphi _{1,0,1,2} , \varphi _{0,1,2,0} ,\varphi _{0,1,2,1} ,\varphi _{0,1,0,2} ,\varphi _{0,1,1,2} ,\varphi _{0,0,2,2} ,\varphi _{1,1,2,0} ,
\]
\[
\left.\varphi _{1,1,2,1} ,\varphi _{1,1,0,2} ,\varphi _{1,1,1,2} ,\varphi _{0,1,2,2} ,\varphi _{1,0,2,2} \right)
\]
из пространства
\[
E_{p}^{\left(1,1,2,2\right)} \equiv L_{p} \left(G\right)\times R\times R\times R\times R\times L_{p} \left(G_{1} \right)\times L_{p} \left(G_{1} \right)\times L_{p} \left(G_{1} \right)\times L_{p} \left(G_{1} \right)\times
\]
\[
\times L_{p} \left(G_{2} \right)\times L_{p} \left(G_{2} \right)\times L_{p} \left(G_{2} \right)\times L_{p} \left(G_{2} \right)\times L_{p} \left(G_{3} \right)\times L_{p} \left(G_{3} \right)\times L_{p} \left(G_{4} \right)\times L_{p} \left(G_{4} \right)\times
\]
\[
\times L_{p} \left(G_{1} \times G_{2} \right)\times L_{p} \left(G_{1} \times G_{2} \right)\times L_{p} \left(G_{1} \times G_{2} \right)\times L_{p} \left(G_{1} \times G_{2} \right)\times L_{p} \left(G_{1} \times G_{3} \right)\times
\]
\[
\times L_{p} \left(G_{1} \times G_{3} \right)\times L_{p} \left(G_{1} \times G_{4} \right)\times L_{p} \left(G_{1} \times G_{4} \right)\times L_{p} \left(G_{2} \times G_{3} \right)\times
\]
\[
\times L_{p} \left(G_{2} \times G_{3} \right)\times L_p \left(G_{2} \times G_{4} \right)\times L_{p} \left(G_{2} \times G_{4} \right)\times L_{p} \left(G_{3} \times G_{4} \right)\times L_{p} \left(G_{1} \times G_{2} \times G_{3} \right)\times
\]
\[
\times L_{p} \left(G_{1} \times G_{2} \times G_{3} \right)\times
%L_{p} \left(G_{1} \times G_{2} \times G_{3} \right)\times
L_{p} \left(G_{1} \times G_{2} \times G_{4} \right)\times
\]
\[
\times L_{p} \left(G_{1} \times G_{2} \times G_{4} \right)\times L_{p} \left(G_{2} \times G_{3} \times G_{4} \right)\times L_{p} \left(G_{1} \times G_{3} \times G_{4} \right).
\]

Заметим, что в пространстве $E_{p}^{\left(1,1,2,2\right)} $ норму будем определять естественным образом, при помощи равенства
\[
\left\| \varphi \right\| _{E_{p}^{\left(1,1,2,2\right)} } =\left\| \varphi _{1,1,2,2} \right\| _{L_{p} \left(G\right)} +\sum _{i_{3} =0}^{1}  \sum _{i_{4} =0}^{1}\left\| \varphi _{0,0,i_{3} ,i_{4} } \right\| _{R}   +\sum _{i_{3} =0}^{1}  \sum _{i_{4} =0}^{1}\left\| \varphi _{1,0,i_{3} ,i_{4} } \right\| _{L_{p} \left(G_{1} \right)}   +
\]
\[
+\sum _{i_{3} =0}^{1}  \sum _{i_{4} =0}^{1}\left\| \varphi _{0,1,i_{3} ,i_{4} } \right\| _{L_{p} \left(G_{2} \right)}   +  \sum _{i_{4} =0}^{1}\left\| \varphi _{0,0,2,i_{4} } \right\| _{L_{p} \left(G_{3} \right)}  +  \sum _{i_{3} =0}^{1}\left\| \varphi _{0,0,i_{3} ,2} \right\| _{L_{p} \left(G_{4} \right)}  +
\]
\[
+\sum _{i_{3} =0}^{1}  \sum _{i_{4} =0}^{1}\left\| \varphi _{1,1,i_{3} ,i_{4} } \right\| _{L_{p} \left(G_{1} \times G_{2} \right)}   +  \sum _{i_{4} =0}^{1}\left\| \varphi _{1,0,2,i_{4} } \right\| _{L_{p} \left(G_{1} \times G_{3} \right)}  +  \sum _{i_{3} =0}^{1}\left\| \varphi _{1,0,i_{3} ,2} \right\| _{L_{p} \left(G_{1} \times G_{4} \right)}  +
\]
\[
+  \sum _{i_{4} =0}^{1}\left\| \varphi _{0,1,2,i_{4} } \right\| _{L_{p} \left(G_{2} \times G_{3} \right)}
 +  \sum _{i_{3} =0}^{1}\left\| \varphi _{0,1,i_{3} ,2} \right\| _{L_{p} \left(G_{2} \times G_{4} \right)}  +
%\sum _{i_{3} =0}^{1}\left\| \varphi _{0,1,i_{3} ,2} \right\| _{L_{p} \left(G_{2} \times G_{4} \right)}  +
\]
\[
+  \left\| \varphi _{0,0,2,2} \right\| _{L_{p} \left(G_{3} \times G_{4} \right)} +  \sum _{i_{4} =0}^{1}\left\| \varphi _{1,1,2,i_{4} } \right\| _{L_{p} \left(G_{1} \times G_{2} \times G_{3} \right)}  +\sum _{i_{3} =0}^{1}\left\| \varphi _{1,1,i_{3} ,2} \right\| _{L_{p} \left(G_{1} \times G_{2} \times G_{4} \right)}  +
\]
\[
+  \left\| \varphi _{0,1,2,2} \right\| _{L_{p} \left(G_{2} \times G_{3} \times G_{4} \right)}   +  \left\| \varphi _{1,0,2,2} \right\| _{L_{p} \left(G_{1} \times G_{3} \times G_{4} \right)}   .
\]

Задачу (1), (4) будем исследовать при помощи интегральных представлений специального вида для функций
$u\left(x\right)\in W_{p}^{\left(1,1,2,2\right)} \left(G\right)$ в виде [19]
%
%Прежде всего отметим, что интегральные представления функций из пространств типа $W_{p}^{\left(1,1,2,2\right)} \left(G\right)$ (из соболевских пространств с доминирующими смешанными производными общего вида) изучены в работах Т.И.Аманова, С.М.Никольского, П.И.Лизор-кина и С.М.Никольского, О.В.Бесова, В.П.Ильина и С.М.Николського, А.Дж.Джабра-илова, С.С.Ахиева, А.М.Наджафова, И.Г.Мамедова и др. Но из них мы будем использовать интегральное представление из, по которому любая функция  единственным образом представима в виде
%
%
\[
u\left(x\right)=\left(Qb\right)  \left(x_{1} ,x_{2} ,x_{3} ,x_{4} \right)\equiv \sum _{i_{3} =0}^{1}  \sum _{i_{4} =0}^{1}x_{3}^{i_{3} } x_{4}^{i_{4} } b_{0,0,i_{3} ,i_{4} }   +
\]
\[
+\sum _{i_{3} =0}^{1}  \sum _{i_{4} =0}^{1}x_{3}^{i_{3} } x_{4}^{i_{4} } \int _{0}^{x_{1} }b_{1,0,i_{3} ,i_{4} } \left(\tau _{1} \right)d\tau _{1}    +\sum _{i_{3} =0}^{1}  \sum _{i_{4} =0}^{1}x_{3}^{i_{3} } x_{4}^{i_{4} } \int _{0}^{x_{2} }b_{0,1,i_{3} ,i_{4} } \left(\tau _{2} \right)d\tau _{2}    +
\]
\[
+  \sum _{i_{4} =0}^{1}x_{4}^{i_{4} } \int _{0}^{x_{3} }\left(x_{3} -\tau _{3} \right)b_{0,0,2,i_{4} } \left(\tau _{3} \right)d\tau _{3}   +\sum _{i_{3} =0}^{1}x_{3}^{i_{3} } \int _{0}^{x_{4} }\left(x_{4} -\tau _{4} \right)b_{0,0,i_{3} ,2} \left(\tau _{4} \right)d\tau _{4}   +
\]
\[
+\sum _{i_{3} =0}^{1}  \sum _{i_{4} =0}^{1}x_{3}^{i_{3} } x_{4}^{i_{4} } \int _{0}^{x_{1} }\int _{0}^{x_{2} }b_{1,1,i_{3} ,i_{4} } \left(\tau _{1} ,\tau _{2} \right)d\tau _{1} d\tau _{2}   +
  \]\[
  +  \sum _{i_{4} =0}^{1}x_{4}^{i_{4} } \int _{0}^{x_{1} }\int _{0}^{x_{3} }\left(x_{3} -\tau _{3} \right)b_{1,0,2,i_{4} } \left(\tau _{1} ,\tau _{3} \right)d\tau _{1} d\tau _{3}    +
\]
\[
+  \sum _{i_{3} =0}^{1}x_{3}^{i_{3} } \int _{0}^{x_{1} }\int _{0}^{x_{4} }\left(x_{4} -\tau _{4} \right)b_{1,0,i_{3} ,2} \left(\tau _{1} ,\tau _{4} \right)d\tau _{1} d\tau _{4}   +\]\[ +  \sum _{i_{4} =0}^{1}x_{4}^{i_{4} } \int _{0}^{x_{2} }\int _{0}^{x_{3} }\left(x_{3} -\tau _{3} \right)b_{0,1,2,i_{4} } \left(\tau _{2} ,\tau _{3} \right)d\tau _{2} d\tau _{3}    +
\]
\[
+  \sum _{i_{3} =0}^{1}x_{3}^{i_{3} } \int _{0}^{x_{2} }\int _{0}^{x_{4} }\left(x_{4} -\tau _{4} \right)b_{0,1,i_{3} ,2} \left(\tau _{2} ,\tau _{4} \right)d\tau _{2} d\tau _{4}  +\]\[  +  \int _{0}^{x_{3} }\int _{0}^{x_{4} }\left(x_{3} -\tau _{3} \right)\left(x_{4} -\tau _{4} \right)b_{0,0,2,2} \left(\tau _{3} ,\tau _{4} \right)d\tau _{3} d\tau _{4}   +
\]
\[
+  \sum _{i_{4} =0}^{1}x_{4}^{i_{4} } \int _{0}^{x_{1} }\int _{0}^{x_{2} }\int _{0}^{x_{3} }\left(x_{3} -\tau _{3} \right)b_{1,1,2,i_{4} } \left(\tau _{1} ,\tau _{2} ,\tau _{3} \right)d\tau _{1} d\tau _{2} d\tau _{3}     +
\]
\[
+  \sum _{i_{3} =0}^{1}x_{3}^{i_{3} } \int _{0}^{x_{1} }\int _{0}^{x_{2} }\int _{0}^{x_{4} }\left(x_{4} -\tau _{4} \right)b_{1,1,i_{3} ,2} \left(\tau _{1} ,\tau _{2} ,\tau _{4} \right)d\tau _{1} d\tau _{2} d\tau _{4}     +
\]
\[
+\int _{0}^{x_{2} }\int _{0}^{x_{3} }\int _{0}^{x_{4} }\left(x_{3} -\tau _{3} \right)\left(x_{4} -\tau _{4} \right)b_{0,1,2,2} \left(\tau _{2} ,\tau _{3} ,\tau _{4} \right)d\tau _{2} d\tau _{3} d\tau _{4}    +
\]
\[
+\int _{0}^{x_{1} }\int _{0}^{x_{3} }\int _{0}^{x_{4} }\left(x_{3} -\tau _{3} \right)\left(x_{4} -\tau _{4} \right)b_{1,0,2,2} \left(\tau _{1} ,\tau _{3} ,\tau _{4} \right)d\tau _{1} d\tau _{3} d\tau _{4}    +
\]
\begin{equation} \label{GrindEQ__6_}
+\int _{0}^{x_{1} }\int _{0}^{x_{2} }\int _{0}^{x_{3} }\int _{0}^{x_{4} }\left(x_{3} -\tau _{3} \right)\left(x_{4} -\tau _{4} \right)b_{1,1,2,2} \left(\tau _{1} ,\tau _{2} ,\tau _{3} ,\tau _{4} \right)d\tau _{1} d\tau _{2} d\tau _{3} d\tau _{4}     ,
\end{equation}
посредством единственного элемента
\[
b=\left(b_{1,1,2,2} ,b_{0,0,0,0} ,b_{0,0,1,0} ,b_{0,0,0,1} ,b_{0,0,1,1} ,b_{1,0,0,0} ,b_{1,0,1,0} ,b_{1,0,0,1} ,b_{1,0,1,1} ,b_{0,1,0,0} ,b_{0,1,1,0} ,b_{0,1,0,1} ,\right.
\]
\[
b_{0,1,1,1} ,b_{0,0,2,0} ,b_{0,0,2,1} ,b_{0,0,0,2} ,b_{0,0,1,2} ,b_{1,1,0,0} ,b_{1,1,1,0} ,b_{1,1,0,1} ,b_{1,1,1,1} ,b_{1,0,2,0} ,b_{1,0,2,1} ,b_{1,0,0,2} ,b_{1,0,1,2} ,
\]
\[
\left. b_{0,1,2,0} ,b_{0,1,2,1} ,b_{0,1,0,2} ,b_{0,1,1,2} ,b_{0,0,2,2} ,b_{1,1,2,0} ,b_{1,1,2,1} ,b_{1,1,0,2} ,b_{1,1,1,2} ,b_{0,1,2,2} ,b_{1,0,2,2} \right)\in E_{p}^{\left(1,1,2,2\right)} .
\]

При этом существуют положительные постоянные $M_{1}^{0} $ и $M_{2}^{0} $ такие, что
\begin{equation} \label{GrindEQ__7_}
M_{1}^{0} \left\| b\right\| _{E_{p}^{\left(1,1,2,2\right)} } \le \left\|   \left(Qb\right)  \left(x_{1} ,x_{2} ,x_{3} ,x_{4} \right)  \right\| _{W_{p}^{\left(1,1,2,2\right)} \left(G\right)} \le M_{2}^{0} \left\| b\right\| _{E_{p}^{\left(1,1,2,2\right)} } ,
\end{equation}
для любой $b\in E_{p}^{\left(1,1,2,2\right)} $.

Очевидно, что оператор $Q:  E_{p}^{\left(1,1,2,2\right)} \to W_{p}^{\left(1,1,2,2\right)} \left(G\right)$ является линейным ограниченным оператором. Неравенство (7) показывает, что оператор $Q$ имеет также ограниченный обратный оператор определенную на пространстве $W_{p}^{\left(1,1,2,2\right)} \left(G\right)$. Следовательно, оператор $Q$ задает гомеоморфизм между банаховыми пространствами $E_{p}^{\left(1,1,2,2\right)} $ и $W_{p}^{\left(1,1,2,2\right)} \left(G\right)$. Поэтому решение уравнения (5) эквивалентно решению уравнения
\begin{equation} \label{GrindEQ__8_}
VQb=\varphi  .
\end{equation}

Уравнение (8) будем называть каноническим видом уравнения (5).

Кроме того, формула (6) показывает, что любая функция $u\in W_{p}^{\left(1,1,2,2\right)} \left(G\right)$ имеет следы:
\[
u\left(0,0,0,0\right), D_{3} u\left(0,0,0,0\right),  D_{4} u\left(0,0,0,0\right),  D_{3} D_{4} u\left(0,0,0,0\right),  D_{1} u\left(x_{1} ,0,0,0\right),
\]
\[
D_{1} D_{3} u\left(x_{1} ,0,0,0\right),D_{1} D_{4} u\left(x_{1} ,0,0,0\right),  D_{1} D_{3} D_{4} u\left(x_{1} ,0,0,0\right),  D_{2} u\left(0,x_{2} ,0,0\right),
\]
\[
D_{2} D_{3} u\left(0,x_{2} ,0,0\right),  D_{2} D_{4} u\left(0,x_{2} ,0,0\right),D_{2} D_{3} D_{4} u\left(0,x_{2} ,0,0\right),  D_{3}^{2} u\left(0,0,x_{3} ,0\right),
\]
\[
 D_{3}^{2} D_{4} u\left(0,0,x_{3} ,0\right),  D_{4}^{2} u\left(0,0,0,x_{4} \right),  D_{3} D_{4}^{2} u\left(0,0,0,x_{4} \right),D_{1} D_{2} u\left(x_{1} ,x_{2} ,0,0\right),
 \]\[
 D_{1} D_{2} D_{3} u\left(x_{1} ,x_{2} ,0,0\right),  D_{1} D_{2} D_{4} u\left(x_{1} ,x_{2} ,0,0\right),  D_{1} D_{2} D_{3} D_{4} u\left(x_{1} ,x_{2} ,0,0\right),
\]
\[
D_{1} D_{3}^{2} u\left(x_{1} ,0,x_{3} ,0\right),  D_{1} D_{3}^{2} D_{4} u\left(x_{1} ,0,x_{3} ,0\right),  D_{1} D_{4}^{2} u\left(x_{1} ,0,0,x_{4} \right),
\]
\[
D_{1} D_{3} D_{4}^{2} u\left(x_{1} ,0,0,x_{4} \right),D_{2} D_{3}^{2} u\left(0,x_{2} ,x_{3} ,0\right),  D_{2} D_{3}^{2} D_{4} u\left(0,x_{2} ,x_{3} ,0\right),
\]
\[
 D_{2} D_{4}^{2} u\left(0,x_{2} ,0,x_{4} \right),  D_{2} D_{3} D_{4}^{2} u\left(0,x_{2} ,0,x_{4} \right),D_{3}^{2} D_{4}^{2} u\left(0,0,x_{3} ,x_{4} \right),
 \]\[
 D_{1} D_{2} D_{3}^{2} u\left(x_{1} ,x_{2} ,x_{3} ,0\right),  D_{1} D_{2} D_{3}^{2} D_{4} u\left(x_{1} ,x_{2} ,x_{3} ,0\right),  D_{1} D_{2} D_{4}^{2} u\left(x_{1} ,x_{2} ,0,x_{4} \right),
\]
\[
D_{1} D_{2} D_{3} D_{4}^{2} u\left(x_{1} ,x_{2} ,0,x_{4} \right),  D_{2} D_{3}^{2} D_{4}^{2} u\left(0,x_{2} ,x_{3} ,x_{4} \right),  D_{1} D_{3}^{2} D_{4}^{2} u\left(x_{1} ,0,x_{3} ,x_{4} \right)
\]
и операции взятия этих следов непрерывны из $W_{p}^{\left(1,1,2,2\right)} \left(G\right)$ в
\[
R, R, R, R, L_{p} \left(G_{1} \right), L_{p} \left(G_{1} \right), L_{p} \left(G_{1} \right), L_{p} \left(G_{1} \right), L_{p} \left(G_{2} \right),
L_{p} \left(G_{2} \right), L_{p} \left(G_{2} \right), L_{p} \left(G_{2} \right),
\]\[
 L_{p} \left(G_{3} \right), L_{p} \left(G_{3} \right), L_{p} \left(G_{4} \right), L_{p} \left(G_{4} \right),
L_{p} \left(G_{1}\times G_{2} \right), L_{p} \left(G_{1}\times G_{2} \right), L_{p} \left(G_{1}\times G_{2} \right),
\]\[
L_{p} \left(G_{1}\times G_{2} \right), L_{p} \left(G_{1}\times G_{3} \right),
L_{p} \left(G_{1}\times G_{3} \right), L_{p} \left(G_{1}\times G_{4} \right), L_{p} \left(G_{1}\times G_{4} \right),
\]\[
L_{p} \left(G_{2}\times G_{3} \right), L_{p} \left(G_{2}\times G_{3} \right),
L_p\left(G_{2}\times G_{4} \right), L_{p} \left(G_{2}\times G_{4} \right), L_{p} \left(G_{3}\times G_{4} \right),
\]\[
L_{p} \left(G_{1}\times G_{2}\times G_{3} \right),
L_{p} \left(G_{1} \times G_{2} \times G_{3} \right), L_{p} \left(G_{1} \times G_{2} \times G_{4} \right),
\]\[
  L_{p} \left(G_{1} \times G_{2} \times G_{4} \right), L_{p} \left(G_{2} \times G_{3} \times G_{4} \right), L_{p} \left(G_{1} \times G_{3} \times G_{4} \right)
\]
соответственно. Далее для этих следов справедливы также равенства:
\[
u\left(0,0,0,0\right)=b_{0,0,0,0} ;  D_{3} u\left(0,0,0,0\right)=b_{0,0,1,0} ;  D_{4} u\left(0,0,0,0\right)=b_{0,0,0,1} ;
\]
\[
D_{3} D_{4} u\left(0,0,0,0\right)=b_{0,0,1,1} ;  D_{1} u\left(x_{1} ,0,0,0\right)=b_{1,0,0,0} \left(x_{1} \right);  D_{1} D_{3} u\left(x_{1} ,0,0,0\right)=b_{1,0,1,0} \left(x_{1} \right);
\]
\[
D_{1} D_{4} u\left(x_{1} ,0,0,0\right)=b_{1,0,0,1} \left(x_{1} \right);  D_{1} D_{3} D_{4} u\left(x_{1} ,0,0,0\right)=b_{1,0,1,1} \left(x_{1} \right);
\]\[
D_{2} u\left(0,x_{2} ,0,0\right)=b_{0,1,0,0} \left(x_{2} \right);  D_{2} D_{3} u\left(0,x_{2} ,0,0\right)=b_{0,1,1,0} \left(x_{2} \right);  \]
\[
D_{2} D_{4} u\left(0,x_{2} ,0,0\right)=b_{0,1,0,1} \left(x_{2} \right);
D_{2} D_{3} D_{4} u\left(0,x_{2} ,0,0\right)=b_{0,1,1,1} \left(x_{2} \right);
\]\[
D_{3}^{2} u\left(0,0,x_{3} ,0\right)=b_{0,0,2,0} \left(x_{3} \right);  D_{3}^{2} D_{4} u\left(0,0,x_{3} ,0\right)=b_{0,0,2,1} \left(x_{3} \right);
\]\[
D_{4}^{2} u\left(0,0,0,x_{4} \right)=b_{0,0,0,2} \left(x_{4} \right);  D_{3} D_{4}^{2} u\left(0,0,0,x_{4} \right)=b_{0,0,1,2} \left(x_{4} \right);
\]
\[
D_{1} D_{2} u\left(x_{1} ,x_{2} ,0,0\right)=b_{1,1,0,0} \left(x_{1} ,x_{2} \right);  D_{1} D_{2} D_{3} u\left(x_{1} ,x_{2} ,0,0\right)=b_{1,1,1,0} \left(x_{1} ,x_{2} \right);
\]\[
D_{1} D_{2} D_{4} u\left(x_{1} ,x_{2} ,0,0\right)=b_{1,1,0,1} \left(x_{1} ,x_{2} \right);  D_{1} D_{2} D_{3} D_{4} u\left(x_{1} ,x_{2} ,0,0\right)=b_{1,1,1,1} \left(x_{1} ,x_{2} \right);
\]
\[
D_{1} D_{3}^{2} u\left(x_{1} ,0,x_{3} ,0\right)=b_{1,0,2,0} \left(x_{1} ,x_{3} \right);  D_{1} D_{3}^{2} D_{4} u\left(x_{1} ,0,x_{3} ,0\right)=b_{1,0,2,1} \left(x_{1} ,x_{3} \right);
\]\[
D_{1} D_{4}^{2} u\left(x_{1} ,0,0,x_{4} \right)=b_{1,0,0,2} \left(x_{1} ,x_{4} \right);  D_{1} D_{3} D_{4}^{2} u\left(x_{1} ,0,0,x_{4} \right)=b_{1,0,1,2} \left(x_{1} ,x_{4} \right);
\]
\[
D_{2} D_{3}^{2} u\left(0,x_{2} ,x_{3} ,0\right)=b_{0,1,2,0} \left(x_{2} ,x_{3} \right);  D_{2} D_{3}^{2} D_{4} u\left(0,x_{2} ,x_{3} ,0\right)=b_{0,1,2,1} \left(x_{2} ,x_{3} \right);
\]\[
D_{2} D_{4}^{2} u\left(0,x_{2} ,0,x_{4} \right)=b_{0,1,0,2} \left(x_{2} ,x_{4} \right);  D_{2} D_{3} D_{4}^{2} u\left(0,x_{2} ,0,x_{4} \right)=b_{0,1,1,2} \left(x_{2} ,x_{4} \right);
\]
\[
D_{3}^{2} D_{4}^{2} u\left(0,0,x_{3} ,x_{4} \right)=b_{0,0,2,2} \left(x_{3} ,x_{4} \right);  D_{1} D_{2} D_{3}^{2} u\left(x_{1} ,x_{2} ,x_{3} ,0\right)=b_{1,1,2,0} \left(x_{1} ,x_{2} ,x_{3} \right);
\]\[
D_{1} D_{2} D_{3}^{2} D_{4} u\left(x_{1} ,x_{2} ,x_{3} ,0\right)=b_{1,1,2,1} \left(x_{1} ,x_{2} ,x_{3} \right);  D_{1} D_{2} D_{4}^{2} u\left(x_{1} ,x_{2} ,0,x_{4} \right)=b_{1,1,0,2} \left(x_{1} ,x_{2} ,x_{4} \right);
\]
\[
D_{1} D_{2} D_{3} D_{4}^{2} u\left(x_{1} ,x_{2} ,0,x_{4} \right)=b_{1,1,1,2} \left(x_{1} ,x_{2} ,x_{4} \right);  D_{2} D_{3}^{2} D_{4}^{2} u\left(0,x_{2} ,x_{3} ,x_{4} \right)=b_{0,1,2,2} \left(x_{2} ,x_{3} ,x_{4} \right);
  \]
  \[
  D_{1} D_{3}^{2} D_{4}^{2} u\left(x_{1} ,0,x_{3} ,x_{4} \right)=b_{1,0,2,2} \left(x_{1} ,x_{3} ,x_{4} \right).
\]
\

\centerline{\bf 3. ЭКВИВАЛЕНТНОЕ  ИНТЕГРАЛЬНОЕ УРАВНЕНИЕ}

Задачу (1), (4) мы будем изучать при помощи интегрального представления (6) функций $u\in W_{p}^{\left(1,1,2,2\right)} \left(G\right)$. Формула (6) показывает, что функция $u\in W_{p}^{\left(1,1,2,2\right)} \left(G\right)$, удовлетворяющая условиям (4), имеет вид:
\[
u\left(x\right)=g_{0} \left(x\right)+\iiiint\limits_{G} R_{0} \left(\tau _{1} ,\tau _{2} ,\tau _{3} ,\tau _{4} ;  x_{1} ,x_{2} ,x_{3} ,x_{4} \right)D_{1} D_{2} D_{3}^{2} D_{4}^{2} u\left(\tau _{1} ,\tau _{2} ,\tau _{3} ,\tau _{4} \right)d\tau _{1} d\tau _{2} d\tau _{3} d\tau _{4}   ,
\]
где
\[
g_{0} \left(x\right)=\sum _{i_{3} =0}^{1}  \sum _{i_{4} =0}^{1}x_{3}^{i_{3} } x_{4}^{i_{4} } \varphi _{0,0,i_{3} ,i_{4} }   +
\sum _{i_{3} =0}^{1}  \sum _{i_{4} =0}^{1}x_{3}^{i_{3} } x_{4}^{i_{4} } \int _{0}^{x_{1} }\varphi _{1,0,i_{3} ,i_{4} } \left(\tau _{1} \right)d\tau _{1}+
\]
\[
   +\sum _{i_{3} =0}^{1}  \sum _{i_{4} =0}^{1}x_{3}^{i_{3} } x_{4}^{i_{4} } \int _{0}^{x_{2} }\varphi _{0,1,i_{3},i_{4} } \left(\tau _{2} \right)d\tau _{2}    +
\sum _{i_{4} =0}^{1}x_{4}^{i_{4} } \int _{0}^{x_{3} }\left(x_{3} -\tau _{3} \right)\varphi _{0,0,2,i_{4} } \left(\tau _{3} \right)d\tau _{3}  +
\]
\[
 +\sum _{i_{3} =0}^{1}x_{3}^{i_{3} } \int _{0}^{x_{4} }\left(x_{4} -\tau _{4} \right)\varphi _{0,0,i_{3} ,2} \left(\tau _{4} \right)d\tau _{4}   +
\sum _{i_{3} =0}^{1}  \sum _{i_{4} =0}^{1}x_{3}^{i_{3} } x_{4}^{i_{4} } \int _{0}^{x_{1} }\int _{0}^{x_{2} }\varphi _{1,1,i_{3} ,i_{4} } \left(\tau _{1} ,\tau _{2} \right)d\tau _{1} d\tau _{2}     +\]
\[
+ \sum _{i_{4} =0}^{1}x_{4}^{i_{4} } \int _{0}^{x_{1} }\int _{0}^{x_{3} }\left(x_{3} -\tau _{3} \right)\varphi _{1,0,2,i_{4} } \left(\tau _{1} ,\tau _{3} \right)d\tau _{1} d\tau _{3}    +
\]
\[
+  \sum _{i_{3} =0}^{1}x_{3}^{i_{3} } \int _{0}^{x_{1} }\int _{0}^{x_{4} }\left(x_{4} -\tau _{4} \right)\varphi _{1,0,i_{3} ,2} \left(\tau _{1} ,\tau _{4} \right)d\tau _{1} d\tau _{4}    +\]
\[
+ \sum _{i_{4} =0}^{1}x_{4}^{i_{4} } \int _{0}^{x_{2} }\int _{0}^{x_{3} }\left(x_{3} -\tau _{3} \right)\varphi _{0,1,2,i_{4} } \left(\tau _{2} ,\tau _{3} \right)d\tau _{2} d\tau _{3}    +
\]
\[
+  \sum _{i_{3} =0}^{1}x_{3}^{i_{3} } \int _{0}^{x_{2} }\int _{0}^{x_{4} }\left(x_{4} -\tau _{4} \right)\varphi _{0,1,i_{3} ,2} \left(\tau _{2} ,\tau _{4} \right)d\tau _{2} d\tau _{4}    +\]
\[
+ \int _{0}^{x_{3} }\int _{0}^{x_{4} }\left(x_{3} -\tau _{3} \right)\left(x_{4} -\tau _{4} \right)\varphi _{0,0,2,2} \left(\tau _{3} ,\tau _{4} \right)d\tau _{3} d\tau _{4}   +
\]
\[
+  \sum _{i_{4} =0}^{1}x_{4}^{i_{4} } \int _{0}^{x_{1} }\int _{0}^{x_{2} }\int _{0}^{x_{3} }\left(x_{3} -\tau _{3} \right)\varphi _{1,1,2,i_{4} } \left(\tau _{1} ,\tau _{2} ,\tau _{3} \right)d\tau _{1} d\tau _{2} d\tau _{3}     +
\]
\[
+  \sum _{i_{3} =0}^{1}x_{3}^{i_{3} } \int _{0}^{x_{1} }\int _{0}^{x_{2} }\int _{0}^{x_{4} }\left(x_{4} -\tau _{4} \right)\varphi _{1,1,i_{3} ,2} \left(\tau _{1} ,\tau _{2} ,\tau _{4} \right)d\tau _{1} d\tau _{2} d\tau _{4}     +
\]
\[
+\int _{0}^{x_{2} }\int _{0}^{x_{3} }\int _{0}^{x_{4} }\left(x_{3} -\tau _{3} \right)\left(x_{4} -\tau _{4} \right)\varphi _{0,1,2,2} \left(\tau _{2} ,\tau _{3} ,\tau _{4} \right)d\tau _{2} d\tau _{3} d\tau _{4}    +
\]
\[
+\int _{0}^{x_{1} }\int _{0}^{x_{3} }\int _{0}^{x_{4} }\left(x_{3} -\tau _{3} \right)\left(x_{4} -\tau _{4} \right)\varphi _{1,0,2,2} \left(\tau _{1} ,\tau _{3} ,\tau _{4} \right)d\tau _{1} d\tau _{3} d\tau _{4}
\]
и
\[
R_{0} \left(\tau _{1} ,\tau _{2} ,\tau _{3} ,\tau _{4} ;  x_{1} ,x_{2} ,x_{3} ,x_{4} \right)=\left(x_{3} -\tau _{3} \right)\left(x_{4} -\tau _{4} \right)\prod _{i=1}^{4}\theta \left(x_{i} -\tau _{i} \right) ,
\]
причем $\theta \left(\xi \right)$ является функцией Хевисайда на $R$, т.е.
\[
\theta \left(\xi \right)=\left\{\begin{array}{l} {1,\quad \xi >0} \\ {0,\quad \xi \le 0.} \end{array}\right.
\]

Тогда после замены $u=g_{0} +\hat{u}$, где
\[
\hat{u}\left(x\right)=\iiiint\limits _{G} R_{0} \left(\tau _{1} ,\tau _{2} ,\tau _{3} ,\tau _{4} ;  x_{1} ,x_{2} ,x_{3} ,x_{4} \right)D_{1} D_{2} D_{3}^{2} D_{4}^{2} u\left(\tau _{1} ,\tau _{2} ,\tau _{3} ,\tau _{4} \right)d\tau _{1} d\tau _{2} d\tau _{3} d\tau _{4}
\]
уравнение (1) можно записать в виде
\begin{equation} \label{GrindEQ__9_}
\left(V_{1,1,2,2} \hat{u}\right)\left(x\right)=\hat{Z}\left(x\right),
\end{equation}
где
\[
\hat{Z}\left(x\right)=\varphi _{1,1,2,2} -V_{1,1,2,2} g_{0} .
\]

Производные функции $\hat{u}$ можно вычислить посредством равенств
\[
D_{1} \hat{u}\left(x\right)=\int _{0}^{x_{2} }\int _{0}^{x_{3} }\int _{0}^{x_{4} }\left(x_{3} -\tau _{3} \right)\left(x_{4} -\tau _{4} \right)D_{1} D_{2} D_{3}^{2} D_{4}^{2} u\left(x_{1} ,\tau _{2} ,\tau _{3} ,\tau _{4} \right)d\tau _{2} d\tau _{3} d\tau _{4}    ,
\]
\[
D_{1} D_{3} \hat{u}\left(x\right)=\int _{0}^{x_{2} }\int _{0}^{x_{3} }\int _{0}^{x_{4} }\left(x_{4} -\tau _{4} \right)D_{1} D_{2} D_{3}^{2} D_{4}^{2} u\left(x_{1} ,\tau _{2} ,\tau _{3} ,\tau _{4} \right)d\tau _{2} d\tau _{3} d\tau _{4}    ,
\]
\[
D_{1} D_{4} \hat{u}\left(x\right)=\int _{0}^{x_{2} }\int _{0}^{x_{3} }\int _{0}^{x_{4} }\left(x_{3} -\tau _{3} \right)D_{1} D_{2} D_{3}^{2} D_{4}^{2} u\left(x_{1} ,\tau _{2} ,\tau _{3} ,\tau _{4} \right)d\tau _{2} d\tau _{3} d\tau _{4}    ,
\]
\[
D_{1} D_{3} D_{4} \hat{u}\left(x\right)=\int _{0}^{x_{2} }\int _{0}^{x_{3} }\int _{0}^{x_{4} }D_{1} D_{2} D_{3}^{2} D_{4}^{2} u\left(x_{1} ,\tau _{2} ,\tau _{3} ,\tau _{4} \right)d\tau _{2} d\tau _{3} d\tau _{4}    ,
\]
\[
D_{2} \hat{u}\left(x\right)=\int _{0}^{x_{1} }\int _{0}^{x_{3} }\int _{0}^{x_{4} }\left(x_{3} -\tau _{3} \right)\left(x_{4} -\tau _{4} \right)D_{1} D_{2} D_{3}^{2} D_{4}^{2} u\left(\tau _{1} ,x_{2} ,\tau _{3} ,\tau _{4} \right)d\tau _{1} d\tau _{3} d\tau _{4}    ,
\]
\[
D_{2} D_{3} \hat{u}\left(x\right)=\int _{0}^{x_{1} }\int _{0}^{x_{3} }\int _{0}^{x_{4} }\left(x_{4} -\tau _{4} \right)D_{1} D_{2} D_{3}^{2} D_{4}^{2} u\left(\tau _{1} ,x_{2} ,\tau _{3} ,\tau _{4} \right)d\tau _{1} d\tau _{3} d\tau _{4}    ,
\]
\[
D_{2} D_{4} \hat{u}\left(x\right)=\int _{0}^{x_{1} }\int _{0}^{x_{3} }\int _{0}^{x_{4} }\left(x_{3} -\tau _{3} \right)D_{1} D_{2} D_{3}^{2} D_{4}^{2} u\left(\tau _{1} ,x_{2} ,\tau _{3} ,\tau _{4} \right)d\tau _{1} d\tau _{3} d\tau _{4}    ,
\]
\[
D_{2} D_{3} D_{4} \hat{u}\left(x\right)=\int _{0}^{x_{1} }\int _{0}^{x_{3} }\int _{0}^{x_{4} }D_{1} D_{2} D_{3}^{2} D_{4}^{2} u\left(\tau _{1} ,x_{2} ,\tau _{3} ,\tau _{4} \right)d\tau _{1} d\tau _{3} d\tau _{4}    ,
\]
\[
D_{3}^{2} \hat{u}\left(x\right)=\int _{0}^{x_{1} }\int _{0}^{x_{2} }\int _{0}^{x_{4} }\left(x_{4} -\tau _{4} \right)D_{1} D_{2} D_{3}^{2} D_{4}^{2} u\left(\tau _{1} ,\tau _{2} ,x_{3} ,\tau _{4} \right)d\tau _{1} d\tau _{2} d\tau _{4}    ,
\]
\[
D_{3}^{2} D_{4} \hat{u}\left(x\right)=\int _{0}^{x_{1} }\int _{0}^{x_{2} }\int _{0}^{x_{4} }D_{1} D_{2} D_{3}^{2} D_{4}^{2} u\left(\tau _{1} ,\tau _{2} ,x_{3} ,\tau _{4} \right)d\tau _{1} d\tau _{2} d\tau _{4}    ,
\]
\[
D_{4}^{2} \hat{u}\left(x\right)=\int _{0}^{x_{1} }\int _{0}^{x_{2} }\int _{0}^{x_{3} }\left(x_{3} -\tau _{3} \right)D_{1} D_{2} D_{3}^{2} D_{4}^{2} u\left(\tau _{1} ,\tau _{2} ,\tau _{3} ,x_{4} \right)d\tau _{1} d\tau _{2} d\tau _{3}    ,
\]
\[
D_{3} D_{4}^{2} \hat{u}\left(x\right)=\int _{0}^{x_{1} }\int _{0}^{x_{2} }\int _{0}^{x_{3} }D_{1} D_{2} D_{3}^{2} D_{4}^{2} u\left(\tau _{1} ,\tau _{2} ,\tau _{3} ,x_{4} \right)d\tau _{1} d\tau _{2} d\tau _{3}    ,
\]
\[
D_{1} D_{2} \hat{u}\left(x\right)=\int _{0}^{x_{3} }\int _{0}^{x_{4} }\left(x_{3} -\tau _{3} \right)\left(x_{4} -\tau _{4} \right)D_{1} D_{2} D_{3}^{2} D_{4}^{2} u\left(x_{1} ,x_{2} ,\tau _{3} ,\tau _{4} \right)d\tau _{3} d\tau _{4}   ,
\]
\[
D_{1} D_{2} D_{3} \hat{u}\left(x\right)=\int _{0}^{x_{3} }\int _{0}^{x_{4} }\left(x_{4} -\tau _{4} \right)D_{1} D_{2} D_{3}^{2} D_{4}^{2} u\left(x_{1} ,x_{2} ,\tau _{3} ,\tau _{4} \right)d\tau _{3} d\tau _{4}   ,
\]
\[
D_{1} D_{2} D_{4} \hat{u}\left(x\right)=\int _{0}^{x_{3} }\int _{0}^{x_{4} }\left(x_{3} -\tau _{3} \right)D_{1} D_{2} D_{3}^{2} D_{4}^{2} u\left(x_{1} ,x_{2} ,\tau _{3} ,\tau _{4} \right)d\tau _{3} d\tau _{4}   ,
\]
\[
D_{1} D_{2} D_{3} D_{4} \hat{u}\left(x\right)=\int _{0}^{x_{3} }\int _{0}^{x_{4} }D_{1} D_{2} D_{3}^{2} D_{4}^{2} u\left(x_{1} ,x_{2} ,\tau _{3} ,\tau _{4} \right)d\tau _{3} d\tau _{4}   ,
\]
\[
D_{1} D_{3}^{2} \hat{u}\left(x\right)=\int _{0}^{x_{2} }\int _{0}^{x_{4} }\left(x_{4} -\tau _{4} \right)D_{1} D_{2} D_{3}^{2} D_{4}^{2} u\left(x_{1} ,\tau _{2} ,x_{3} ,\tau _{4} \right)d\tau _{2} d\tau _{4}   ,
\]
\[
D_{1} D_{3}^{2} D_{4} \hat{u}\left(x\right)=\int _{0}^{x_{2} }\int _{0}^{x_{4} }D_{1} D_{2} D_{3}^{2} D_{4}^{2} u\left(x_{1} ,\tau _{2} ,x_{3} ,\tau _{4} \right)d\tau _{2} d\tau _{4}   ,
\]
\[
D_{1} D_{4}^{2} \hat{u}\left(x\right)=\int _{0}^{x_{2} }\int _{0}^{x_{3} }\left(x_{3} -\tau _{3} \right)D_{1} D_{2} D_{3}^{2} D_{4}^{2} u\left(x_{1} ,\tau _{2} ,\tau _{3} ,x_{4} \right)d\tau _{2} d\tau _{3}   ,
\]
\[
D_{1} D_{3} D_{4}^{2} \hat{u}\left(x\right)=\int _{0}^{x_{2} }\int _{0}^{x_{3} }D_{1} D_{2} D_{3}^{2} D_{4}^{2} u\left(x_{1} ,\tau _{2} ,\tau _{3} ,x_{4} \right)d\tau _{2} d\tau _{3}   ,
\]
\[
D_{2} D_{3}^{2} \hat{u}\left(x\right)=\int _{0}^{x_{1} }\int _{0}^{x_{4} }\left(x_{4} -\tau _{4} \right)D_{1} D_{2} D_{3}^{2} D_{4}^{2} u\left(\tau _{1} ,x_{2} ,x_{3} ,\tau _{4} \right)d\tau _{1} d\tau _{4}   ,
\]
\[
D_{2} D_{3}^{2} D_{4} \hat{u}\left(x\right)=\int _{0}^{x_{1} }\int _{0}^{x_{4} }D_{1} D_{2} D_{3}^{2} D_{4}^{2} u\left(\tau _{1} ,x_{2} ,x_{3} ,\tau _{4} \right)d\tau _{1} d\tau _{4}   ,
\]
\[
D_{2} D_{4}^{2} \hat{u}\left(x\right)=\int _{0}^{x_{1} }\int _{0}^{x_{3} }\left(x_{3} -\tau _{3} \right)D_{1} D_{2} D_{3}^{2} D_{4}^{2} u\left(\tau _{1} ,x_{2} ,\tau _{3} ,x_{4} \right)d\tau _{1} d\tau _{3}   ,
\]
\[
D_{2} D_{3} D_{4}^{2} \hat{u}\left(x\right)=\int _{0}^{x_{1} }\int _{0}^{x_{3} }D_{1} D_{2} D_{3}^{2} D_{4}^{2} u\left(\tau _{1} ,x_{2} ,\tau _{3} ,x_{4} \right)d\tau _{1} d\tau _{3}   ,
\]
\[
D_{3} \hat{u}\left(x\right)=\int _{0}^{x_{1} }\int _{0}^{x_{2} }\int _{0}^{x_{3} }\int _{0}^{x_{4} }\left(x_{4} -\tau _{4} \right)D_{1} D_{2} D_{3}^{2} D_{4}^{2} u\left(\tau _{1} ,\tau _{2} ,\tau _{3} ,\tau _{4} \right)d\tau _{1} d\tau _{2} d\tau _{3} d\tau _{4}     ,
\]
\[
D_{4} \hat{u}\left(x\right)=\int _{0}^{x_{1} }\int _{0}^{x_{2} }\int _{0}^{x_{3} }\int _{0}^{x_{4} }\left(x_{3} -\tau _{3} \right)D_{1} D_{2} D_{3}^{2} D_{4}^{2} u\left(\tau _{1} ,\tau _{2} ,\tau _{3} ,\tau _{4} \right)d\tau _{1} d\tau _{2} d\tau _{3} d\tau _{4}     ,
\]
\[
D_{3} D_{4} \hat{u}\left(x\right)=\int _{0}^{x_{1} }\int _{0}^{x_{2} }\int _{0}^{x_{3} }\int _{0}^{x_{4} }D_{1} D_{2} D_{3}^{2} D_{4}^{2} u\left(\tau _{1} ,\tau _{2} ,\tau _{3} ,\tau _{4} \right)d\tau _{1} d\tau _{2} d\tau _{3} d\tau _{4}     ,
\]
\[
D_{3}^{2} D_{4}^{2} \hat{u}\left(x\right)=\int _{0}^{x_{1} }\int _{0}^{x_{2} }D_{1} D_{2} D_{3}^{2} D_{4}^{2} u\left(\tau _{1} ,\tau _{2} ,x_{3} ,x_{4} \right)d\tau _{1} d\tau _{2}   ,
\]
\[
D_{1} D_{2} D_{3}^{2} \hat{u}\left(x\right)=\int _{0}^{x_{4} }\left(x_{4} -\tau _{4} \right)D_{1} D_{2} D_{3}^{2} D_{4}^{2} u\left(x_{1} ,x_{2} ,x_{3} ,\tau _{4} \right)d\tau _{4}  ,
\]
\[
D_{1} D_{2} D_{3}^{2} D_{4} \hat{u}\left(x\right)=\int _{0}^{x_{4} }D_{1} D_{2} D_{3}^{2} D_{4}^{2} u\left(x_{1} ,x_{2} ,x_{3} ,\tau _{4} \right)d\tau _{4}  ,
\]
\[
D_{1} D_{2} D_{4}^{2} \hat{u}\left(x\right)=\int _{0}^{x_{3} }\left(x_{3} -\tau _{3} \right)D_{1} D_{2} D_{3}^{2} D_{4}^{2} u\left(x_{1} ,x_{2} ,\tau _{3} ,x_{4} \right)d\tau _{3}  ,
\]
\[
D_{1} D_{2} D_{3} D_{4}^{2} \hat{u}\left(x\right)=\int _{0}^{x_{3} }D_{1} D_{2} D_{3}^{2} D_{4}^{2} u\left(x_{1} ,x_{2} ,\tau _{3} ,x_{4} \right)d\tau _{3}  ,
\]
\[
D_{2} D_{3}^{2} D_{4}^{2} \hat{u}\left(x\right)=\int _{0}^{x_{1} }D_{1} D_{2} D_{3}^{2} D_{4}^{2} u\left(\tau _{1} ,x_{2} ,x_{3} ,x_{4} \right)d\tau _{1}  ,
\]
\[
D_{1} D_{3}^{2} D_{4}^{2} \hat{u}\left(x\right)=\int _{0}^{x_{2} }D_{1} D_{2} D_{3}^{2} D_{4}^{2} u\left(x_{1} ,\tau _{2} ,x_{3} ,x_{4} \right)d\tau _{2}  ,
\]
\[
D_{1} D_{2} D_{3}^{2} D_{4}^{2} \hat{u}\left(x\right)=D_{1} D_{2} D_{3}^{2} D_{4}^{2} u\left(x_{1} ,x_{2} ,x_{3} ,x_{4} \right).
\]

Теперь доминирующую производную рассмотрим как неизвестную функцию, иначе говоря, произведем замену
\[
D_{1} D_{2} D_{3}^{2} D_{4}^{2} u\left(x_{1} ,x_{2} ,x_{3} ,x_{4} \right)=b\left(x_{1} ,x_{2} ,x_{3} ,x_{4} \right).
\]

Тогда уравнение (9) можно записать в виде:
\[
\left(Nb\right)\left(x\right)\equiv b\left(x_{1} ,x_{2} ,x_{3} ,x_{4} \right)+
\int _{0}^{x_{1} }\int _{0}^{x_{2} }\int _{0}^{x_{3} }\int _{0}^{x_{4} }a_{0,0,0,0} \left(x_{1} ,x_{2} ,x_{3} ,x_{4} \right)
\times
\]
\[
\times R_{0} \left(\tau _{1} ,\tau _{2} ,\tau _{3} ,\tau _{4} ;  x_{1} ,x_{2} ,x_{3} ,x_{4} \right)b\left(\tau _{1} ,\tau _{2} ,\tau _{3} ,\tau _{4} \right)d\tau _{1} d\tau _{2} d\tau _{3} d\tau _{4}     +
\]
\[
+\int _{0}^{x_{1} }\int _{0}^{x_{2} }\int _{0}^{x_{3} }\int _{0}^{x_{4} }a_{0,0,1,0} \left(x_{1} ,x_{2} ,x_{3} ,x_{4} \right)\left(x_{4} -\tau _{4} \right)b\left(\tau _{1} ,\tau _{2} ,\tau _{3} ,\tau _{4} \right)d\tau _{1} d\tau _{2} d\tau _{3} d\tau _{4}     +
\]
\[
+\int _{0}^{x_{1} }\int _{0}^{x_{2} }\int _{0}^{x_{3} }\int _{0}^{x_{4} }a_{0,0,0,1} \left(x_{1} ,x_{2} ,x_{3} ,x_{4} \right)\left(x_{3} -\tau _{3} \right)b\left(\tau _{1} ,\tau _{2} ,\tau _{3} ,\tau _{4} \right)d\tau _{1} d\tau _{2} d\tau _{3} d\tau _{4}     +
\]
\[
+\int _{0}^{x_{1} }\int _{0}^{x_{2} }\int _{0}^{x_{3} }\int _{0}^{x_{4} }a_{0,0,1,1} \left(x_{1} ,x_{2} ,x_{3} ,x_{4} \right)b\left(\tau _{1} ,\tau _{2} ,\tau _{3} ,\tau _{4} \right)d\tau _{1} d\tau _{2} d\tau _{3} d\tau _{4}     +
\]
\[
+\int _{0}^{x_{2} }\int _{0}^{x_{3} }\int _{0}^{x_{4} }a_{1,0,0,0} \left(x_{1} ,x_{2} ,x_{3} ,x_{4} \right)\left(x_{3} -\tau _{3} \right)\left(x_{4} -\tau _{4} \right)b\left(x_{1} ,\tau _{2} ,\tau _{3} ,\tau _{4} \right)d\tau _{2} d\tau _{3} d\tau _{4}    +
\]
\[
+\int _{0}^{x_{2} }\int _{0}^{x_{3} }\int _{0}^{x_{4} }a_{1,0,1,0} \left(x_{1} ,x_{2} ,x_{3} ,x_{4} \right)\left(x_{4} -\tau _{4} \right)b\left(x_{1} ,\tau _{2} ,\tau _{3} ,\tau _{4} \right)d\tau _{2} d\tau _{3} d\tau _{4}    +
\]
\[
+\int _{0}^{x_{2} }\int _{0}^{x_{3} }\int _{0}^{x_{4} }a_{1,0,0,1} \left(x_{1} ,x_{2} ,x_{3} ,x_{4} \right)\left(x_{3} -\tau _{3} \right)b\left(x_{1} ,\tau _{2} ,\tau _{3} ,\tau _{4} \right)d\tau _{2} d\tau _{3} d\tau _{4}    +
\]
\[
+\int _{0}^{x_{2} }\int _{0}^{x_{3} }\int _{0}^{x_{4} }a_{1,0,1,1} \left(x_{1} ,x_{2} ,x_{3} ,x_{4} \right)b\left(x_{1} ,\tau _{2} ,\tau _{3} ,\tau _{4} \right)d\tau _{2} d\tau _{3} d\tau _{4}    +
\]
\[
+\int _{0}^{x_{1} }\int _{0}^{x_{3} }\int _{0}^{x_{4} }a_{0,1,0,0} \left(x_{1} ,x_{2} ,x_{3} ,x_{4} \right)\left(x_{3} -\tau _{3} \right)\left(x_{4} -\tau _{4} \right)b\left(\tau _{1} ,x_{2} ,\tau _{3} ,\tau _{4} \right)d\tau _{1} d\tau _{3} d\tau _{4}    +
\]
\[
+\int _{0}^{x_{1} }\int _{0}^{x_{3} }\int _{0}^{x_{4} }a_{0,1,1,0} \left(x_{1} ,x_{2} ,x_{3} ,x_{4} \right)\left(x_{4} -\tau _{4} \right)b\left(\tau _{1} ,x_{2} ,\tau _{3} ,\tau _{4} \right)d\tau _{1} d\tau _{3} d\tau _{4}    +
\]
\[
+\int _{0}^{x_{1} }\int _{0}^{x_{3} }\int _{0}^{x_{4} }a_{0,1,0,1} \left(x_{1} ,x_{2} ,x_{3} ,x_{4} \right)\left(x_{3} -\tau _{3} \right)b\left(\tau _{1} ,x_{2} ,\tau _{3} ,\tau _{4} \right)d\tau _{1} d\tau _{3} d\tau _{4}    +
\]
\[
+\int _{0}^{x_{1} }\int _{0}^{x_{3} }\int _{0}^{x_{4} }a_{0,1,1,1} \left(x_{1} ,x_{2} ,x_{3} ,x_{4} \right)b\left(\tau _{1} ,x_{2} ,\tau _{3} ,\tau _{4} \right)d\tau _{1} d\tau _{3} d\tau _{4}    +
\]
\[
+\int _{0}^{x_{1} }\int _{0}^{x_{2} }\int _{0}^{x_{4} }a_{0,0,2,0} \left(x_{1} ,x_{2} ,x_{3} ,x_{4} \right)\left(x_{4} -\tau _{4} \right)b\left(\tau _{1} ,\tau _{2} ,x_{3} ,\tau _{4} \right)d\tau _{1} d\tau _{2} d\tau _{4}    +
\]
%\[
%+\int _{0}^{x_{1} }\int _{0}^{x_{2} }\int _{0}^{x_{4} }a_{0,0,2,0} \left(x_{1} ,x_{2} ,x_{3} ,x_{4} \right)\left(x_{4} -\tau _{4} \right)b\left(\tau _{1} ,\tau _{2} ,x_{3} ,\tau _{4} \right)d\tau _{1} d\tau _{2} d\tau _{4}    +
%\]
\[
+\int _{0}^{x_{1} }\int _{0}^{x_{2} }\int _{0}^{x_{4} }a_{0,0,2,1} \left(x_{1} ,x_{2} ,x_{3} ,x_{4} \right)b\left(\tau _{1} ,\tau _{2} ,x_{3} ,\tau _{4} \right)d\tau _{1} d\tau _{2} d\tau _{4}    +
\]
\[
+\int _{0}^{x_{1} }\int _{0}^{x_{2} }\int _{0}^{x_{3} }a_{0,0,0,2} \left(x_{1} ,x_{2} ,x_{3} ,x_{4} \right)\left(x_{3} -\tau _{3} \right)b\left(\tau _{1} ,\tau _{2} ,\tau _{3} ,x_{4} \right)d\tau _{1} d\tau _{2} d\tau _{3}    +
\]
\[
+\int _{0}^{x_{1} }\int _{0}^{x_{2} }\int _{0}^{x_{3} }a_{0,0,1,2} \left(x_{1} ,x_{2} ,x_{3} ,x_{4} \right)b\left(\tau _{1} ,\tau _{2} ,\tau _{3} ,x_{4} \right)d\tau _{1} d\tau _{2} d\tau _{3}    +
\]
\[
+\int _{0}^{x_{3} }\int _{0}^{x_{4} }a_{1,1,0,0} \left(x_{1} ,x_{2} ,x_{3} ,x_{4} \right)\left(x_{3} -\tau _{3} \right)\left(x_{4} -\tau _{4} \right)b\left(x_{1} ,x_{2} ,\tau _{3} ,\tau _{4} \right)d\tau _{3} d\tau _{4}   +
\]
\[
+\int _{0}^{x_{3} }\int _{0}^{x_{4} }a_{1,1,1,0} \left(x_{1} ,x_{2} ,x_{3} ,x_{4} \right)\left(x_{4} -\tau _{4} \right)b\left(x_{1} ,x_{2} ,\tau _{3} ,\tau _{4} \right)d\tau _{3} d\tau _{4}   +
\]
\[
+\int _{0}^{x_{3} }\int _{0}^{x_{4} }a_{1,1,0,1} \left(x_{1} ,x_{2} ,x_{3} ,x_{4} \right)\left(x_{3} -\tau _{3} \right)b\left(x_{1} ,x_{2} ,\tau _{3} ,\tau _{4} \right)d\tau _{3} d\tau _{4}   +
\]
\[
+\int _{0}^{x_{3} }\int _{0}^{x_{4} }a_{1,1,1,1} \left(x_{1} ,x_{2} ,x_{3} ,x_{4} \right)b\left(x_{1} ,x_{2} ,\tau _{3} ,\tau _{4} \right)d\tau _{3} d\tau _{4}   +
\]
\[
+\int _{0}^{x_{2} }\int _{0}^{x_{4} }a_{1,0,2,0} \left(x_{1} ,x_{2} ,x_{3} ,x_{4} \right)\left(x_{4} -\tau _{4} \right)b\left(x_{1} ,\tau _{2} ,x_{3} ,\tau _{4} \right)d\tau _{2} d\tau _{4}   +
\]
\[
+\int _{0}^{x_{2} }\int _{0}^{x_{4} }a_{1,0,2,1} \left(x_{1} ,x_{2} ,x_{3} ,x_{4} \right)b\left(x_{1} ,\tau _{2} ,x_{3} ,\tau _{4} \right)d\tau _{2} d\tau _{4}   +
\]
\[
+\int _{0}^{x_{2} }\int _{0}^{x_{3} }a_{1,0,0,2} \left(x_{1} ,x_{2} ,x_{3} ,x_{4} \right)\left(x_{3} -\tau _{3} \right)b\left(x_{1} ,\tau _{2} ,\tau _{3} ,x_{4} \right)d\tau _{2} d\tau _{3}   +
\]
\[
+\int _{0}^{x_{2} }\int _{0}^{x_{3} }a_{1,0,1,2} \left(x_{1} ,x_{2} ,x_{3} ,x_{4} \right)b\left(x_{1} ,\tau _{2} ,\tau _{3} ,x_{4} \right)d\tau _{2} d\tau _{3}   +
\]
\[
+\int _{0}^{x_{1} }\int _{0}^{x_{4} }a_{0,1,2,0} \left(x_{1} ,x_{2} ,x_{3} ,x_{4} \right)\left(x_{4} -\tau _{4} \right)b\left(\tau _{1} ,x_{2} ,x_{3} ,\tau _{4} \right)d\tau _{1} d\tau _{4}   +
\]
\[
+\int _{0}^{x_{1} }\int _{0}^{x_{4} }a_{0,1,2,1} \left(x_{1} ,x_{2} ,x_{3} ,x_{4} \right)b\left(\tau _{1} ,x_{2} ,x_{3} ,\tau _{4} \right)d\tau _{1} d\tau _{4}   +
\]
\[
+\int _{0}^{x_{1} }\int _{0}^{x_{3} }a_{0,1,0,2} \left(x_{1} ,x_{2} ,x_{3} ,x_{4} \right)\left(x_{3} -\tau _{3} \right)b\left(\tau _{1} ,x_{2} ,\tau _{3} ,x_{4} \right)d\tau _{1} d\tau _{3}   +
\]
\[
+\int _{0}^{x_{1} }\int _{0}^{x_{3} }a_{0,1,1,2} \left(x_{1} ,x_{2} ,x_{3} ,x_{4} \right)b\left(\tau _{1} ,x_{2} ,\tau _{3} ,x_{4} \right)d\tau _{1} d\tau _{3}   +
\]
\[
+\int _{0}^{x_{1} }\int _{0}^{x_{2} }a_{0,0,2,2} \left(x_{1} ,x_{2} ,x_{3} ,x_{4} \right)b\left(\tau _{1} ,\tau _{2} ,x_{3} ,x_{4} \right)d\tau _{1} d\tau _{2}   +
\]
\[
+\int _{0}^{x_{4} }a_{1,1,2,0} \left(x_{1} ,x_{2} ,x_{3} ,x_{4} \right)\left(x_{4} -\tau _{4} \right)b\left(x_{1} ,x_{2} ,x_{3} ,\tau _{4} \right)d\tau _{4}  +
\]
\[
+\int _{0}^{x_{4} }a_{1,1,2,1} \left(x_{1} ,x_{2} ,x_{3} ,x_{4} \right)b\left(x_{1} ,x_{2} ,x_{3} ,\tau _{4} \right)d\tau _{4}  +
\]
\[
+\int_0^{x_3}a_{1,1,0,2}\left(x_1,x_2,x_3,x_4\right)\left(x_3-\tau_3\right)b\left(x_1,x_2,\tau_3,x_4\right)d\tau_3+
\]
\[
+\int _{0}^{x_{3} }a_{1,1,1,2} \left(x_{1} ,x_{2} ,x_{3} ,x_{4} \right)b\left(x_{1} ,x_{2} ,\tau _{3} ,x_{4} \right)d\tau _{3}  +
\]
\[
+\int _{0}^{x_{1} }a_{0,1,2,2} \left(x_{1} ,x_{2} ,x_{3} ,x_{4} \right)b\left(\tau _{1} ,x_{2} ,x_{3} ,x_{4} \right)d\tau _{1}  +
\]
\begin{equation} \label{GrindEQ__10_}
+\int _{0}^{x_{2} }a_{1,0,2,2} \left(x_{1} ,x_{2} ,x_{3} ,x_{4} \right)b\left(x_{1} ,\tau _{2} ,x_{3} ,x_{4} \right)d\tau _{2}  =\hat{Z}\left(x\right),  \ \     x\in G
\end{equation}

Оператор $N$ уравнения (10) линеен. Используя условия наложенные на коэффициенты $a_{i_{1} ,i_{2} ,i_{3} ,i_{4} } $, можно доказать, что этот оператор является ограниченным оператором из $L_{p} \left(G\right)$ в $L_{p} \left(G\right),    1\le p\le \infty $.

\textbf{Определение.} {\it Если задача (1), (4) для любого
\[
\varphi =\left(\varphi _{1,1,2,2} ,\varphi _{0,0,0,0} ,\varphi _{0,0,1,0} ,\varphi _{0,0,0,1} ,\varphi _{0,0,1,1} ,\varphi _{1,0,0,0} ,\varphi _{1,0,1,0} ,\varphi _{1,0,0,1} ,\varphi _{1,0,1,1} ,\varphi _{0,1,0,0} ,\varphi _{0,1,1,0} ,\right.
\]
\[
\varphi _{0,1,0,1} ,\varphi _{0,1,1,1} ,\varphi _{0,0,2,0} ,\varphi _{0,0,2,1} ,\varphi _{0,0,0,2} ,\varphi _{0,0,1,2} ,\varphi _{1,1,0,0} ,\varphi _{1,1,1,0} ,\varphi _{1,1,0,1} ,\varphi _{1,1,1,1} ,\varphi _{1,0,2,0} ,\varphi _{1,0,2,1} ,
\]
\[
\varphi _{1,0,0,2} ,\varphi _{1,0,1,2} ,\varphi _{0,1,2,0} ,\varphi _{0,1,2,1} ,\varphi _{0,1,0,2} ,\varphi _{0,1,1,2} ,\varphi _{0,0,2,2} ,\varphi _{1,1,2,0} ,\varphi _{1,1,2,1} ,
\]\[
\left. \varphi _{1,1,0,2} ,\varphi _{1,1,1,2} ,\varphi _{0,1,2,2} ,\varphi _{1,0,2,2} \right)\in E_{p}^{\left(1,1,2,2\right)}
\]
имеет единственное решение $u\in W_{p}^{\left(1,1,2,2\right)} \left(G\right)$ такое, что
\[
\left\| u\right\| _{W_{p}^{\left(1,1,2,2\right)} \left(G\right)} \le M_{1} \left\| \varphi\right\| _{E_{p}^{\left(1,1,2,2\right)} } ,
\]
то будем говорить, что оператор $V$ задачи (1), (4) (или уравнения (5)) задает гомеоморфизм из $W_{p}^{\left(1,1,2,2\right)} \left(G\right)$ на $E_{p}^{\left(1,1,2,2\right)} $  или задача (1), (4) везде корректно разрешима. Здесь $M_{1} $ постоянная не зависящее от $\varphi $.}

Очевидно, что, если оператор $V$ задачи (1), (4) задает гомеоморфизм из $W_{p}^{\left(1,1,2,2\right)} \left(G\right)$ на $E_{p}^{\left(1,1,2,2\right)} $, то существует ограниченный обратный оператор $V^{-1} :\linebreak E_{p}^{\left(1,1,2,2\right)} \to W_{p}^{\left(1,1,2,2\right)} \left(G\right)$.

В современной теории дифференциальных уравнений особое значение имеет вопрос о выявлении классов задач, операторы которых осуществляют гомеоморфизм между определенными парами банаховых пространств. Такие гомеоморфизмы выявлены в работах Ю.М.Березанского и Я.А.Ройтберга [20], Н.В.Житарашу [21], \linebreak С.С.Ахиева [22], И.Г.Мамедова [23] и др. для некоторых классов дифференциальных уравнений с частными производными.

Оператор $N$ является вольтерровым оператором относительно точки $\left(0,0,0,0\right)$. Это означает, что если функции $b_{1} ,  b_{2} \in L_{p} \left(G\right)$ в области
\[
G_{\left(x_{1} ,x_{2} ,x_{3} ,x_{4} \right)} =\left(0,x_{1} \right)\times \left(0,x_{2} \right)\times \left(0,x_{3} \right)\times \left(0,x_{4} \right)
\]
удовлетворяют условию
\[
b_{1} \left(\tau _{1} ,\tau _{2} ,\tau _{3} ,\tau _{4} \right)=b_{2} \left(\tau _{1} ,\tau _{2} ,\tau _{3} ,\tau _{4} \right),
\]
то выполняется также условие
\[
\left(Nb_{1} \right)  \left(\tau _{1} ,\tau _{2} ,\tau _{3} ,\tau _{4} \right)=\left(Nb_{2} \right)  \left(\tau _{1} ,\tau _{2} ,\tau _{3} ,\tau _{4} \right)
\]
почти для всех $\left(\tau _{1} ,\tau _{2} ,\tau _{3} ,\tau _{4} \right)\in G_{\left(x_{1} ,x_{2} ,x_{3} ,x_{4} \right)} $, где $\left(x_{1} ,x_{2} ,x_{3} ,x_{4} \right)\in G$ произвольная точка.

Используя вольтерровость оператора $N$, при помощи, например, метода последовательных приближений можно доказать, что уравнение (10) для любой правой части $\hat{Z}\in L_{p} \left(G\right)$ имеет единственное решение $b\in L_{p} \left(G\right)$, где $1\le p\le \infty $ и это решение удовлетворяет условию
\[
\left\| b\right\| _{L_{p} \left(G\right)} \le M_{2} \left\| \hat{Z}\right\| _{L_{p} \left(G\right)} ,
\]
где $M_{2} $ постоянное не зависящее от $\hat{Z}$. Далее, очевидно, что если $\varphi _{1,1,2,2} \in L_{p} \left(G\right)$, то $\hat{Z}\in L_{p} \left(G\right)$. Кроме того, если $b\in L_{p} \left(G\right)$ есть решение уравнения (10), то решение задачи (1), (4) можно найти при помощи равенства
\[
u\left(x\right)=g_{0} \left(x\right)+\int _{0}^{x_{1} }\int _{0}^{x_{2} }\int _{0}^{x_{3} }\int _{0}^{x_{4} }b\left(\tau _{1} ,\tau _{2} ,\tau _{3} ,\tau _{4} \right)R_{0} \left(\tau _{1} ,\tau _{2} ,\tau _{3} ,\tau _{4} ;  x_{1} ,x_{2} ,x_{3} ,x_{4} \right)d\tau _{1} d\tau _{2} d\tau _{3} d\tau _{4}     .
\]

Поэтому справедлива

\textbf{Теорема.} {\it Оператор задачи (1), (4) задает гомеоморфизм из $W_{p}^{\left(1,1,2,2\right)} \left(G\right)$ на $E_{p}^{\left(1,1,2,2\right)} $.}\\

\centerline{\bf 4. ВЫВОДЫ}

Постановка задачи (1), (4) обладает рядом преимуществ:

1) в этой постановке не требуется никаких дополнительных условий согласования;

2) именно, такая постановка порождает гомеоморфизм между двумя банаховыми пространствами $W_p^{1,1,2,2}(G)$
и $E_p^{1,1,2,2}$;

3) эту задачу можно рассматривать как задачу сформулированную по следам в пространстве С.Л.Соболева $W_p^{1,1,2,2}(G)$;

4) в этой постановке рассматриваемое уравнение является обобщениям многих модельных уравнений некоторых процессов (например,  уравнения влагопереноса, уравнения теплопроводности, уравнения Аллера, уравнения Буссинеска-Лява, уравнения Манжерона, трехмерного телеграфного уравнения и т.д.).

Отметим, что применяя методику, приводимую в статье [24], доказывается интегральное представление решения задачи (1), (4) используя понятия фундаментального решения.

\newpage

\begin{titlepage}

\centerline{\bf Сведения об авторе}

\

{\bf Мамедов Ильгар Гурбат оглы }

Азербайджан, AZ 1141, г. Баку, ул. Б. Вагабзаде, 9,

Институт Кибернетики НАН Азербайджана, www.cyber.az.

Тел.: (99412) 539 28 26,

Факс.: (99412) 539 28 26,

e-mail: ilgar-mammadov@rambler.ru.

Должность: ведущий научный сотрудник.

Ученая степень: кандидат физико-математических наук.

Уч. звание: доцент.
\end{titlepage}
\end{document}